\documentclass[10pt,twoside,reqno]{amsart} %%%%%%%%%%%%%%%%%% reqno: c'est right equation number
\usepackage{amsmath, amsfonts, amsthm, amssymb}
\usepackage[francais,english]{babel}
\usepackage{caption}

\usepackage[latin1]{inputenc}
\usepackage{amsmath}
\usepackage{amsfonts}
\usepackage{amsthm}
\usepackage[]{graphics}
\usepackage{color}
\usepackage[T1]{fontenc}
 \usepackage{setspace}
\usepackage{algorithm}
\usepackage{algorithmic}

\usepackage[T1]{fontenc}
\usepackage[latin1]{inputenc}
\usepackage{amsmath,latexsym}
\usepackage{here}
\usepackage{graphicx}
\usepackage{multirow}
\usepackage{amssymb}

\definecolor{red}{rgb}{0.82,0.00,0.00}  % rouge

\numberwithin{equation}{section}
\input epsf
 

\newcommand{\ds}{\displaystyle}
\def\R{{\rm I\hspace{-0.50ex}R} }
\def\P{{\rm I\hspace{-0.50ex}P} }

\def\w{{\bf w}}
\def\r{{\bf r}}
\def\f{{\bf f}}
\def\v{{ \bf v}}
\def\V{{ \bf V}}
\def\u{{\bf u}}

\def\W{{\bf W}}

\def\0{{\bf 0}}
\def\p{p}
\def\q{q}

\def\x{{\bf{x}}}
\def\g{{\bf g}}
\def\G{{\bf G}}

\def\n{{\bf n}}
\def\div{\operatorname{div}}
\def\curl{\operatorname{\bf curl}}
\newtheorem{lem}{Lemma}[section]

\newtheorem{thm}[lem]{Theorem}
\newtheorem{rmq}[lem]{Remark}
\newtheorem{Asmp}[lem]{Assumption}

\newtheorem{defi}[lem]{Definition}
\newtheorem{prop}[lem]{Proposition}
\newtheorem{propri}[lem]{Property}
\renewcommand{\epsilon}{\varepsilon}
\def\twoplot[#1]#2#3#4#5{
\begin{figure}[hbt]
\begin{multicols}{2}
\begin{center}
    \includegraphics*[#1]{#2}
    \caption{\label{#2} #4}
\end{center}
\begin{center}
    \includegraphics*[#1]{#3}
    \caption{\label{#3} #5}
\end{center}
\end{multicols}
\end{figure}
}

\begin{document}

\bibliographystyle{plain}

\title[A posteriori error estimates for the BDF problem]{ A posteriori error estimates for the Brinkman-Darcy-Forchheimer problem}
\author[ SAYAH  ]{Toni Sayah$^{\dagger}$}
\thanks{ \today.
\newline
$^{\dagger}$  Laboratoire de "Math\'ematiques et applications", Unit\'e de recherche "Math\'ematoqies et Mod\'elisation", CAR, Facult\'e des sciences, Universit\'e Saint-Joseph, Lebanon.
\newline
toni.sayah@usj.edu.lb.}

%\date{}

\begin{abstract}
In this paper, we study the {\it a posteriori} error estimate corresponding to the Brinkman-Darcy-Forchheimer problem. We introduce the variational formulation  discretised by using the finite element method. Then, we establish an {\it a posteriori} error estimation with two types of error indicators related to the discretization and to the linearization. Finally, numerical investigations are shown and discussed.\\

\noindent $ \mathbf{ \sc Keywords:}$   Brinkman-Darcy-Forchheimer problem, finite element method, {\it a posteriori} error estimation.

\end{abstract}

\maketitle

\section{Introduction}\label{intro}
\noindent Let $\Omega$ be a bounded subset of $\R^d$ ($d = 2, 3$) with Lipschitz continuous boundary $\Gamma=\partial \Omega$.  We consider the Brinkman-Darcy-Forchheimer equation (see for instance \cite{BDF1,BDF2,BDF3} )
\begin{equation}\label{E1}
\ds -\div (Re^{-1} \varepsilon \nabla \u - \varepsilon \u \otimes \u)  + \alpha (\varepsilon)  \u + \beta (\varepsilon) |\u| \u + \varepsilon \nabla p =  \varepsilon \f \;\;\;\;\; \mbox{in}\;\;\; \Omega,
\end{equation}
with the divergence constraint
\begin{equation}\label{E2}
\ds \div (\varepsilon \u) = 0 \;\;\;\;\; \mbox{in}\;\;\; \Omega,
\end{equation}
and the boundary condition
\begin{equation}\label{E3}
\ds \u = \g \;\;\;\;\; \mbox{on}\;\; \Gamma.
\end{equation}
Here $\f$ is an external force field, $\g$ is a given function on $\Gamma$, $\u$ represents the velocity, $p$ represents the pressure, $|.|$ denotes the Euclidean norm, $|\u|^2 = \u \cdot \u$. The positive function $\varepsilon$ represents the porosity of the domain and varies spatially in general. The functions $\alpha$ and $\beta$ represent the Darcy and Forchheimer terms and $Re$ designates the Reynolds number.\\
\noindent We denote by Problem $(P_1)$ the system of equations \eqref{E1}, \eqref{E2}, \eqref{E3}. \\

The Darcy-Brinkman-Forchheimer (DBF) model for porous media is obtained from the incompressible Navier-Stokes
equation in a porous media. The importance of the BDF model is that it can be used to model porous media with relatively large Reynolds numbers, since the Darcy equation is considered to be a suitable for small range of
Reynolds numbers. For the derivation of the equations, their limitations, modelling, and homogenization questions in porous media we refer to \cite{vaf,BDF1, Hornung, Garibotti, Matossian, Winterberg, Diersch, Nield, Whitaker}. \\
For the mathematical study of the steady-state DBF equation, we can refer to \cite{BDFtheor1, Kaloni, Skrzypacz1} for the existence and uniqueness of the solution of the problem with inhomogeneous Dirichlet boundary condition. In \cite{BDFtheornum1}, the authors treated the steady state BDF problem with mixed boundary condition and proved that the considered problem has a unique solution if the source terms are small enough. Then, the convergence of a Taylor-Hood finite element approximation using a finite element interpolation of the porosity is then proved under
similar smallness assumptions. Some optimal error estimates are obtained and some numerical experiments are performed.\\

On the other hand, {\it a posteriori} error estimation consists in bounding the error between the exact and the numerical solutions with a sum of local indicators at each cell of the mesh. In order to obtain a more accurate solution with a low additional CPU cost, we can adapt the mesh with respect to the local values of the indicators. The {\it a posteriori} error estimate is optimal if each indicator can be bounded by the local error in the solution around the corresponding element. The {\it a posteriori} error analysis was first introduced by Babu$\check{s}$ka and Rheinboldt \cite{BaRh} and then further developed, among others, by Verf\"{u}rth \cite{Verfurth} or Ainsworth and Oden \cite{AiOd}. In \cite{Elakkad}, El Akkad, El Khalifi and Guessous proposed a discretization of the incompressible Navier-Stokes equations by mixed finite elements with {\it a posteriori} error estimation of the computed solutions. Other works on {\it a posteriori} estimations for the stationary Navier-Stokes system are proposed in \cite{Jin,john,BDMS,BHV,joanna} and for the non-stationary Navier-Stokes system in \cite{BcSt,NaSt}. \\

In this paper we consider Problem $(P_1)$, and introduce the variational formulation and  some corresponding properties established in \cite{BDFtheor1,BDFtheornum1}. Then, we introduce the discrete variational formulation by using the finite element method and by stabilizing the convection term. An iterative numerical scheme is then introduced and the corresponding convergence is studied. Next, an {\it a posteriori} error estimate is established and finally, numerical investigations are performed in order to show the importance of the adapted mesh method versus the uniform mesh method. \\

This paper is organised as follow :
\begin{itemize}
\item Section 2 describes the problem and the weak formulation.
\item Section 3 is devoted to study the discretization and the convergence of the proposed iterative schemes.
\item Section 4 is devoted the {\it a posteriori} error studies.
\item In Section 5 we show numerical investigations.
\end{itemize}
\section{Notations and weak formulation}
In order to introduce the variational formulation, we  recall  some classical Sobolev spaces and  their  properties.

Let $\alpha=(\alpha_1,\alpha_2, \dots \alpha_d)$ be a $d$-uple of non negative
integers, set $|\alpha|=\ds \sum_{i=1}^d \alpha_i$, and define the partial
derivative $\partial^\alpha$ by
$$
\partial^{\alpha}=\ds \frac{\partial^{|\alpha|}}{\partial x_1^{\alpha_1}\partial
x_2^{\alpha_2}\dots\partial x_d^{\alpha_d}}.
$$
Then, for any positive integer $m$ and number $q\geq 1$, we recall
the classical Sobolev space
\begin{equation}
\label{eq:Wm,p}
W^{m,q}(\Omega)=\{v \in L^q(\Omega);\,\forall\,|\alpha|\leq
m,\;\partial^{\alpha} v \in L^q(\Omega)\},
\end{equation}
equipped with the semi-norm
\begin{equation}
\label{eq:semnormWm,p}
|v|_{W^{m,q}(\Omega)}=\big\{\sum_{|\alpha|=m} \int_{\Omega}
|\partial^{\alpha} v|^{  q} \,d\x\,\big\}^{\frac{1}{q}}
\end{equation}
and the norm
\begin{equation}
\label{eq:normWm,p}
\|v\|_{ W^{m,q}(\Omega)}=\big\{\sum_{0\leq k\leq m}
|v|_{ W^{k,q}(\Omega)}^q \big\}^{ \frac{1}{q}}.
\end{equation}
When $ q=2$, this space is the Hilbert space $H^m(\Omega)$.
In
particular, the scalar product of $L^2(\Omega)$ is denoted by
$(.,.)$.
The definitions of these spaces are extended straightforwardly to
vectors, with the same notation, but with the following
modification for the norms in the non-Hilbert case. Let $\v$ be a
vector valued function {  and we define the norm}
\begin{equation}
\label{eq:normLp} \|\v\|_{ L^q(\Omega)}= \big(\int _{\Omega}
|\v|^{ q}\,d\x\, \big)^{ \frac{1}{q}},
\end{equation}
where $|.|$ denotes the Euclidean vector norm.\\
We have the useful lemma:
\begin{lem}\label{lemmapp}
For any $ p \le 6 $, there exists positive constant $S_p$ and $C_p$ only depending on $\Omega$ such that
\begin{equation}\label{L2Xo1} \forall  \v \in H^1\Omega)^d, \quad
||\v||_{L^p(\Omega)^d} \le S_p ||\v||_{H^1(\Omega)^d}
\end{equation}
and
\begin{equation}\label{L2Xo2} \forall  \v \in H^1_0(\Omega)^d, \quad
||\v||_{L^p(\Omega)^d} \le C_p |\v|_{H^1_0(\Omega)^d}.
\end{equation}
\end{lem}
%
%\noindent Let us now introduce the following technical lemma:
%
%{\lem \label{relatinf} For all $x,y \in \R$ and $q\in \R^{+}$, the following bound holds:
%\[
%(|x|^{q} x - |y|^{q} y)(x-y) \ge 0.
%\]
%}
%
%
\subsection{Variational formulation}
We refer to \cite{BDFtheor1,BDFtheornum1} for all the properties and the details of the weak formulation corresponding to Problem $(P_1)$ presented in this section. We suppose that $\f \in L^2(\Omega)^d$ and we begin by giving some assumptions on the porosity, the Darcy and the Forchheimer functions of the fluid.
{ \Asmp \label{assumpdata}  We assume that the porosity, the Darcy and Forchheimer terms satisfy the next set of assumptions:
\begin{itemize}
\item $\varepsilon \in L^\infty(\Omega) \cap W^{1,r}(\Omega)$ with $r > d$ and $0 < \varepsilon_0 \le \varepsilon_(\x) \le 1$ in $\Omega$.
\item $\alpha$ and $\beta$ are positive in $[\varepsilon_0,1]$ and differentiable in $]\varepsilon_0,1[$.
\item $\alpha'$ and $\beta'$ are positive and bounded in $]\varepsilon_0,1[$.
\end{itemize}
}
Assumption \ref{assumpdata} gives that there exists $\alpha_{m}, \alpha_{M}, \beta_{m}$ and $\beta_{M}$ such that $\alpha \in [\alpha_{m}, \alpha_{M}]$ and $\beta \in [\beta_{m}, \beta_{M}]$.\\

By using the incompressibility condition $\div(\varepsilon \u)=0$, the non-linear term can be written as
\[
\div (\varepsilon \u \otimes \u) = \div(\varepsilon \u) \u + \varepsilon (\u \cdot \nabla ) \u = \varepsilon (\u \cdot \nabla ) \u.
\]
Let us define the spaces
\[
\begin{array}{ll}
\medskip
X = H^1(\Omega)^d,\quad X_0 = H^1_0(\Omega)^d,\\
M= L^2_0(\Omega).
\end{array}
\]
We denote by $X_0^{'}$ the dual space of $X_0$.\\

Problem $(P_1)$ is equivalent to the following variational formulation \cite{BDFtheor1,BDFtheornum1}: Find $(\u,p)\in X \times M$ such that $\u_\Gamma = \g$ and
\begin{equation}\label{V1}
\left\{
\begin{array}{ll}
\medskip
\forall \v \in X_0, \quad a(\u,\v) + c(\u,\u,\v) - b(\v,p) = \ell(\v),\\
\forall q\in M, \quad b(\u,q)=0,
\end{array}
\right.
\end{equation}
where
\[
\begin{array}{rcl}
\medskip
a(\u,\v) &=& \ds \frac{1}{Re} (\varepsilon \nabla \u, \nabla \v) + (\alpha(\varepsilon) \u, \v),\\
\medskip
c(\u,\v,\w) &=& (\varepsilon (\u \cdot \nabla) \v , \w) + (\beta(\varepsilon) |\u| \v, \w),\\
\medskip
b(\u,p) &=& (\div(\varepsilon \u), p),\\
\ell(\v) &=& \langle \varepsilon \f, \v\rangle.
\end{array}
\]
To deal with the inhomogeneous Dirichlet condition, we introduce an extension ${\G}\in X$ of $\g$ such that $\div(\varepsilon {\G}) = 0$ whose existence is provided by \cite{BDFtheornum1} (Lemma 16). From this, $\u = \w + {\G}$ where $\w \in X_0$. Thus Problem \eqref{V1} allows us to deduce that $\w$ is the solution of the following problem:
\begin{equation}\label{W1}
\left\{
\begin{array}{ll}
\medskip
\forall \v \in X_0, \quad a(\w,\v) + c(\w +{\G},\w +{\G},\v)   - b(\v,p) = \ell(\v)  -a({\G},v)) ,\\
\forall q\in M, \quad b(\w,q)=0.
\end{array}
\right.
\end{equation}
With the previous definition of the bilinear form $b$ and if $\varepsilon \in L^\infty(\Omega) \cap W^{1,r(d)}(\Omega)$ with $r(2)>2$ and $r(3)=3$, we have the following inf-sup condition (see \cite{BDFtheornum1}):
\[
\underset{\v \in X_0}{\sup} \; \; \ds \frac{b(\v,q)}{||\v||_{X_0(\Omega)}} \ge \gamma_\varepsilon \, || q ||_{M},
\]
where $\gamma_\varepsilon$ is a constant which depends on $\varepsilon$.\\

We recall the following theorem of the existence and uniqueness of the solution of Problem \eqref{W1}. For the proof, we refer to \cite{BDFtheor1,BDFtheornum1}.
{\thm \label{exuniqtheo} Under Assumption \ref{assumpdata}, Problem \eqref{W1} admits at least one solution $(\u,p)$ in $X \times M$. \\
Furthermore, if the data $\f$ and $\g$ are sufficiently small in the sense that there exists a positive real number $\eta_e$ such that
\[
|| \f ||_{X_0^{'}} + ||\g ||_{H^{1/2}(\Gamma)^d} \le  \eta_e,
\]
then the solution of problem \eqref{W1} is unique and there exists a positive real number $R_e$ such that
\[
||\w||_X + ||p||_M \le R_e.
\]
}
\noindent We note that Theorem \ref{exuniqtheo} gives the existence and uniqueness of a solution to \eqref{V1} since $\u=\w+{\G}$.\\
%where ${\G}$ is the divergence-free lifting of $\g$.\\

%
%
%
For the development of the {\it a posteriori} error estimate, we introduce the Stokes problem defined  as follows: For a given $(\f,t) \in X_0^{'} \times M$, look for $(\bar \w,\bar p)$ in $X_0 \times M$ such that
\begin{equation*}
\label{fVss}
\hspace{-1cm} (S)
\left\{
\begin{array}{lll}
\medskip
\mbox{Find} \ \bar \w \in X_0 \mbox{ and } \bar p \in  M \mbox{ such that:} & &  \\
\medskip
\forall \v \in X_0, \ a (\bar \w, \v)  - b(\v,  \bar p) = \ell(\v), &&  \\
\forall q \in  M, \ \ b(\bar \w, q) = { (t,q)}.
\end{array} \right.
\end{equation*}
We refer to the general abstract framework given by \cite[Chapter I, pargraphe 4]{Girault}
for the existence and the uniqueness of the solution  $(\bar \w, \bar p) \in X_0 \times M$ of Problem $(S)$. \\
\noindent We introduce the following Stokes operator $ \mathcal S$:
 $$\begin{array}{ccccc}
\mathcal S  & : &  X_0^{'} \times M & \to & X_0 \times M \\
 & & (\f,t) & \mapsto & \mathcal S (\f,t)= ( \bar \w,\bar p ),
\end{array}$$

\noindent where $ (\bar \w,\bar p)$ is the solution of  the Stokes problem $(S)$. We have the following bound (see  \cite[Chapter I, paragraphe 4]{Girault} )
\begin{equation}\label{stokes3}
||\mathcal S (\f,t)||_{X_0 \times M} \le \ds \frac{c}{\nu}  ||(\f,t) ||_{X_0^{'}\times M}.
\end{equation}
%%%
\noindent We define also the function $ \mathcal{G} $ given by
$$\begin{array}{ccccc}
\mathcal{G}  & : &  X_0 & \to & X_0^{'}\\
& & \bar \u & \mapsto &  \mathcal{G}(\bar \u)
\end{array}
$$
where for all $\v \in X_0$,
\[
\langle \mathcal{G}(\bar \u), \v \rangle =  \ds a({\G}, \v) + c(\bar \u + {\G}, \bar \u + {\G}, \v) - \ell(\v),
\]
\noindent and we introduce the map $F$ on $ X_0 \times M $ such that for all $ \V=(\v,q)  \in X_0 \times M$, we have
$$ F(\V)= \V+\mathcal S (\mathcal{G}(\v),0).$$
Then, Problem  \eqref{W1} can be equivalently written as
\begin{equation} \label{fu}
F(\W)=\0,
\end{equation}
where $ \W=(\w,p)$ is the solution of \eqref{W1}.
\begin{defi}
We define a non singular solution $ \W$ of Problem \eqref{W1} in the following way:
\begin{enumerate}
\item $ F(\W)=\0.$
\item $ DF(\W)$  is an  isomorphism of $ X_0 \times M$,
\end{enumerate}
\noindent where  $ DF(\W)$  is the Fr\'echet-differential of $F$ in $\W$.
\end{defi}
In the next proposition, we show the Lipschitz property concerning the operator $DF$. For the prove, we refer to \cite[Theorem 11]{BDFtheornum1}.
%
%
%-------------------------------------------
\begin{prop} \label{isom}
There exist a real number $L > 0$ and a neighborhood $\mathcal{V}$ of $\W =(\w,p)$ in $ X_0 \times M $ such that the following Lipschitz property holds: For all  $\V=(\v,\q)$ in $\mathcal{V}$, we have
  \begin{equation}  \label{DFisom}
|| DF(\W)- DF(\V)||_{\mathcal L(X_0 \times M)}  \le  L  || \W-\V||_{X_0 \times M }.
\end{equation}
\end{prop}
\section{Finite element discretization and convergence}
From now on, we assume that $\Omega$ is a polygon when $d=2$ or
polyhedron when $d=3$, so it can be completely meshed.  For
 the space discretization,  we consider a
regular (see Ciarlet~\cite{PGC}) family of triangulations
$( \mathcal{T}_h )_h$ of $\Omega$ which is a set of closed non degenerate
triangles for $d=2$ or tetrahedra for $d=3$, called elements,
satisfying,
\begin{itemize}
\item for each $h$, $\bar{\Omega}$ is the union of all elements of
 $\mathcal{T}_h$;
\item the intersection of two distinct elements of  $\mathcal{T}_h$ is
either empty, a common vertex, or an entire common edge (or face
when $d=3$);
\item the ratio of the diameter   {  $h_{\kappa}$  of an element $\kappa \in \mathcal{T}_h$  to
the diameter  $\rho_\kappa$}  of its inscribed circle when $d=2$ or ball when $d=3$
is bounded by a constant independent of $h$: there exists a positive constant $\sigma$ independent of $h$ such that,
\begin{equation}
\label{eq:reg}
{
\ds \max_{\kappa \in \mathcal{T}_{h}} \frac{h_{ \kappa}}{\rho_{ \kappa}} \le \sigma.}
\end{equation}
\end{itemize}
As usual, $h$ denotes the maximal diameter
of all elements of $\mathcal{T}_{h}$.   To define the finite element functions,  let $r$ be a non negative integer.  For each {  $\kappa$ in $\mathcal{T}_{h}$, we denote by  $\mathbb{P}_r(\kappa)$  the space of restrictions to $\kappa$ of polynomials in $d$ variables and total degree at most  $r$, with a similar notation on the faces or edges of $\kappa$}. For every edge (when $d=2$) or face (when $d=3$) $e$ of the mesh  $\mathcal{T}_h$, we denote by $h_e$ the diameter of $e$.
%\\
%
%\noindent {  We shall use the following inverse   inequality {  \cite{sayahdib}}:  for any dimension $d$,  there exists a constant $C_I$ such that for any polynomial function $v_h$ {  of degree $r$   on $\kappa$},
%for any number $p\geq 2$ and  for any dimension $d$,  there exists a constant $C_I(p)$ such that for any polynomial function $v_h$ {  of degree $r$   on $\kappa$},
%\begin{equation}
%\label{eq:inversin}
% {  \|v_h\|_{L^p(\kappa)} \leq  C_I(p)  h_\kappa^{\frac{d}{p}-\frac{d}{2}}\|v_h\|_{L^2(\kappa)}.}
%\end{equation}
%The constant $C_I(p)$ depends also on the regularity parameter $\sigma$ of \eqref{eq:reg}, but for the sake of simplicity this is not indicated.
%\begin{equation}
%\label{eq:inversin}
%{  \|v_h\|_{L^p(\kappa)} \leq  C_I(p)  h_\kappa^{\frac{d}{p}-\frac{d}{2}}\|v_h\|_{L^2(\kappa)}.}
%\end{equation}
%The constant $C_I(p)$ depends on the regularity parameter $\sigma$ of \eqref{eq:reg}, but for the sake of simplicity this is not indicated.}
%
\subsection{Discretization of the variational problem}
In this section we follow the discretizations introduced in \cite{BDFtheornum1} where the authors introduce a discrete variational formulation corresponding to \eqref{W1} based on the Taylor-Hood finite element \cite{Taulorhood}. But in this work, we use the mini-elements for the discretization.\\

Let $X_{0h} \subset X_0$ and $M_{0h} \subset M$ the discrete spaces corresponding to the velocity and the pressure. In \cite{BDFtheornum1}, the authors consider the Taylor-Hood finite element and introduce and study a discrete variational problem. They also establish the corresponding {\it a priori} error estimate between the exact and approximated solutions. In this work, The choice of the spaces $X_{0h}$ and $M_{0h}$ can be based on the Taylor-Hood finite element (see for instance \cite{BDFtheornum1}) or on the mini-element (see \cite{minielement}). The results established in \cite{BDFtheornum1} concerning the Taylor-Hood element can be easily extended to the mini-element which will be considered in this work.\\

Let $\kappa$ be an element of ${\mathcal T}_h$ with vertices $a_i$, $1\leq i \leq d+1$, and
corresponding barycentric coordinates $\lambda_i$. We denote
by $b_\kappa \in \mathbb{P}_{d+1}(\kappa)$ the basic bubble function
\begin{equation}
\label{eqn:bbl}
b_\kappa(\x)=\lambda_1(\x)...\lambda_{d+1}(\x).
\end{equation}
We observe that $b_\kappa(\x)=0$ on $\partial \kappa$ and that $b_\kappa(\x)>0$
in the interior of $\kappa$.\\
Let $(X_{0h},M_{0h})$ be the discrete spaces defined by
\begin{equation}
\label{eq:Wh2}
X_h = \{\v_h \in (\mathcal{C}^0(\bar{\Omega}))^d;\;\forall\, \kappa \in \mathcal{T}_h,\;
\v_{h}|_\kappa\in {\mathcal{P}(\kappa)}^d\} ,
\end{equation}
\begin{equation}
X_{0h}=X_h \cap X_0,
\end{equation}
\begin{equation}
M_h=\{q_h \in  \mathcal{C}^0(\bar{\Omega});\;\forall\, \kappa \in \mathcal{T}_h,\; q_{h}|_\kappa\in \mathbb{P}_1(\kappa)\},
\end{equation}
\begin{equation}
M_{0h}=M_h \cap L^2_0(\Omega),
\end{equation}
where
$$\mathcal{P}(\kappa)= \mathbb{P}_1(\kappa) \oplus {\rm Vect}\{b_\kappa\},$$
is associated to the discretization in space by the "mini-element" (see \cite{minielement}).\\

Let $\mathcal{V}_{0h}$ be the kernel of the divergence in $X_{0h}$,
\begin{equation}
\label{eq:Vh2} \mathcal{V}_{0h}=\{ \v_h\in
X_{0h};\;\forall q_h \in M_{0h},\ds
\int_{\Omega}  q_h \, \div (\varepsilon_h \v_h)\,d\x=0 \}.
\end{equation}
In order to introduce a discrete scheme associated to \eqref{W1}, I consider some $\varepsilon_h \in M_h$ that approximates $\varepsilon$ in the following sense:
\begin{itemize}
\item $\forall \x \in \Omega, \varepsilon_0 \le \varepsilon_h(\x) \le 1$.
\item $\| \varepsilon - \varepsilon_h \|_{L^\infty (\Omega)} \le C h \| \varepsilon \|_{W^{1,\infty}(\Omega)} $.
\item $\| \nabla \varepsilon - \nabla \varepsilon_h \|_{L^r (\Omega)} \le C h \| \varepsilon \|_{W^{2,r}(\Omega)}, r > d $.
\end{itemize}
Note that these assumptions are satisfied if $\varepsilon \in W^{1,\infty}(\Omega) \cap W^{2,r}(\Omega)$ and if one takes $\varepsilon_h = \mathcal{I}_h \varepsilon$ where $\mathcal{I}_h$ is the global interpolation operator (see \cite{ErnGuerm}, Corollary 1.109 and Corollary 1.110]).\\

The finite element approximation $\u_h$ of the solution to Problem \eqref{V1} is going to be $\u_h = \w_h + \mathcal{I}_{X_h} {\G}$ where $\mathcal{I}_{X_h} : X \rightarrow X_h$ is the finite element interpolate operator
%satisfying the following inequality:
%\begin{equation}\label{stabili}
%\| \mathcal{I}_{X_h} {\G} \|_X \le c_s %\|{\G} \|_X,
%\end{equation}
%where $c_s$ is a constant independent %of $h$. This operator can be for %instance the Cl\'ement operator (see \cite [section IX.3]{BMR} and \cite{clement}).\\
%Thus
and $(\w_h,p_h)\in X_{0h} \times M_{0h}$ is the solution of the following discrete problem associated to \eqref{W1}:
\begin{equation}\label{W1h}
\left\{
\begin{array}{ll}
\medskip
\forall \v_h \in X_{0h}, \quad a_h(\w_h,\v_h) +  c_h(\w_h + {\G},\w_h + {\G},\v_h) -b_h(\v_h,p_h) = \ell_h(\v_h) -a_h({\G},\v_h)) ,\\
\forall q_h\in M_{0h}, \quad b_h(\w_h,q_h)=0,
\end{array}
\right.
\end{equation}
where
\[
\begin{array}{rcl}
\medskip
a_h(\u,\v) &=& \ds \frac{1}{Re} (\varepsilon_h \nabla \u_h, \nabla \v_h) + (\alpha(\varepsilon_h) \u_h, \v_h),\\
\medskip
c_h(\u_h,\v_h,\w_h) &=& (\varepsilon_h (\u_h \cdot \nabla) \v_h , \w_h)  + \ds \frac{1}{2} (\div(\varepsilon_h\, \u_h) \v_h, \w_h)+ (\beta(\varepsilon_h) |\u_h| \v_h, \w_h),\\
\medskip
b_h(\u_h,p_h) &=& (\div(\varepsilon_h \u_h), p_h),\\
\ell_h(\v_h) &=& \langle \varepsilon_h \f, \v_h\rangle.
\end{array}
\]
We introduce the form $d_h$ given by
\[
d_h(\u_h,\v_h,\w_h) = (\varepsilon_h (\u_h \cdot \nabla) \v_h , \w_h)  + \ds \frac{1}{2} (\div(\varepsilon_h\, \u_h) \v_h, \w_h)
\]
and we have
\[
c_h(\u_h,\v_h,\w_h) =d_h(\u_h,\v_h,\w_h) +  (\beta(\varepsilon_h) |\u_h| \v_h, \w_h).
\]
We note that the form $d_h$ is linear with all their variables while the form $c_h$ is linear only with the second and third variables $\v_h$ and $\w_h$. Furthermore, the term $\ds \frac{1}{2} (\div(\varepsilon_h\, \u_h) \v_h, \w_h)$ is called the term of stabilisation since we have for all $\v_h, \w_h \in X_{0h}$,
\[
\ds (\varepsilon_h (\v_h \cdot \nabla) \w_h , \w_h)  + \frac{1}{2} (\div(\varepsilon_h\, \v_h) \w_h, \w_h) = 0
\]
\begin{lem}\label{propch}
For the simplicity of the following calculus, we show some properties concerning the forms $c_h$ and $d_h$. For each $\u_h,\v_h$ and $\w_h$ in $X_{0h}$, we have
%\begin{enumerate}
%\item for each $\u_h,\v_h$ and $\w_h$ in $X_{0h}$, we have,
\begin{equation}\label{rrr1}
\begin{array}{ll}
\medskip
d_h({\G},\w_h,\w_h)=0, \\
\medskip
d_h(\v_h, \w_h,\w_h) = 0,\\
\medskip
c_h(\v_h, \w_h,\w_h) = (\beta(\varepsilon_h) |\u_h| \w_h, \w_h) \ge 0,\\
\medskip
d_h(\u_h, \v_h,\w_h) \le c_d | \u_h |_X | \v_h |_X |\w_h |_X,\\
c_h(\u_h, \v_h,\w_h) \le c_v | \u_h |_X | \v_h |_X |\w_h |_X .
\end{array}
\end{equation}
%
%\item there exists a positive constant $c_g$ such that
%\begin{equation}\label{rrr1g}
%\begin{array}{ll}
%c_h({\G}, \v_h,\w_h) \le c_g \| {\G} %\|_X  | \v_h |_X |\w_h |_X,  \\
%c_h(\v_h,{\G}, \w_h) \le c_g \| {\G} %\|_X  | \v_h |_X |\w_h |_X.
%\end{array}
%\end{equation}
%
%\item there exists a positive constant %$c_{gg}$ such that
%\begin{equation}\label{rrr1gg}
%c_h({\G}, {\G}, \w_h) \le c_{gg} \| %{\G} \|_X^2 |\w_h |_X.
%\end{equation}
%\end{enumerate}
\end{lem}
\noindent \textbf {Proof.} The first two lines of \eqref{rrr1} are straightforward.
The proof of the last three terms are a simple consequence of the definition of the forms $d_h$ and $c_h$, the properties of $\varepsilon_h$ and Lemma \ref{lemmapp} for $p=2,4$.\hfill$\Box$\\
{\rmq We have for each $\u_h,\v_h$ and $\w_h$ in $X_h$,
\begin{equation}\label{rrr1g}
\begin{array}{ll}
%d_h(\v_h, \w_h,\w_h) = 0,\\
%c_h(\v_h, \w_h,\w_h) = (\beta(\varepsilon_h) |\u_h| \w_h, \w_h) \ge 0,\\
\medskip
d_h(\u_h, \v_h,\w_h) \le c_d \| \u_h \|_X \| \v_h \|_X \|\w_h \|_X,\\
c_h(\u_h, \v_h,\w_h) \le c_v \| \u_h \|_X \| \v_h \|_X \|\w_h \|_X.
\end{array}
\end{equation}
}

\noindent We note that the discrete scheme \eqref{W1h} is similar to that introduced in \cite{BDFtheornum1} where the only difference is in the definition of the form $c_h$ which contains here the term of stabilisation $\ds \frac{1}{2} (\div(\varepsilon_h\, \u_h) \v_h, \w_h)$. In fact this term play an important role in the iterative scheme introduced later. The introduction of the stabilisation term change slightly the studies carried out in in \cite{BDFtheornum1}.\\

With the previous definition of $b_h$ we have the following discrete inf-sup condition (see \cite{BDFtheornum1})
\begin{equation}
\label{infsuph1h}
\forall\, q_h\in M_{0h},\; \sup_{\v_h\in X_{0h}} \ds \frac{b_h(\v_h,q_h)}{\|\v_h\|_{X_0}}\geq \gamma^*_\varepsilon  \|q_h\|_{M},
\end{equation}
where $\gamma^*_\varepsilon$ is a positive constant independent of $h$ but depends on $\varepsilon$.\\

In the next theorem, we prove the existence and uniqueness of the solutions of Problem \eqref{W1h}.

\begin{thm}\label{exuniqnum}
Under Assumption \ref{assumpdata}, there exists  positive real number $\eta_n$ (given by \eqref{rrr6}) such that if
\begin{equation}\label{rrr5}
|| {\g}||_{H^{1/2}(\Gamma)^d} < \ds \eta_n,
\end{equation}
Problem \eqref{W1h} admits at least one solution $(\w_h,p_h) \in X_{0h} \times M_{0h}$ such that
\begin{equation}\label{majwh}
|\w_h |_X \le R_1(\f,\g),
\end{equation}
where $R_1(\f,\g)$ is given by \eqref{R1}.\\
Furthermore, if ${\g}$ and $\f$ are sufficiently small such that
\[
\varepsilon_0 \| \f \|_{X_0^{'}}  + c_v || {\g} ||^2_{H^{1/2}(\Gamma)^d}   + \ds (\frac{\varepsilon_0 + 2}{2Re }  + \alpha_M C_2^2)|| {\g} ||_{H^{1/2}(\Gamma)^d}  < \ds \frac{\varepsilon_0^2}{2Re^2(c_d + \beta_M C_2 C_4 S_4)},
\]
then the solution of Problem \eqref{W1h} is unique.
\end{thm}
\noindent \textbf {Proof.}
Let us define the nonlinear mapping $F: \mathcal{V}_{0h} \rightarrow \mathcal{V}_{0h}$ such that for each $\w_h \in \mathcal{V}_{0h}$, we have for all $\v_h \in \mathcal{V}_{0h}$:
\[
\begin{array}{rcl}
\medskip
(\nabla F(\w_h), \nabla \v_h) &=& a_h(\w_h, \v_h) + c_h(\w_h + {\G},\w_h+ {\G},\v_h) - \ell_h(\v_h) + a_h({\G},\v_h).
\end{array}
\]
By using \eqref{rrr1}, we get for each $\v_h \in \mathcal{V}_{0h}$,
\[
\begin{array}{rcl}
\medskip
c_h(\v_h + {\G},\v_h+ {\G},\v_h) &=& d_h(\w_h + {\G},\v_h+ {\G},\v_h)  + (\beta(\varepsilon_h)|\v_h + {\G}| (\v_h+ {\G}),\v_h) \\
\medskip
&=& d_h(\w_h + {\G}, {\G},\v_h)  + (\beta(\varepsilon_h)|\v_h + {\G}| {\G},\v_h)\\
\medskip
&& + (\beta(\varepsilon_h)|\v_h + {\G}| \v_h,\v_h) ,  \\
&=& c_h(\w_h + {\G}, {\G},\v_h)  + (\beta(\varepsilon_h)|\v_h + {\G}| \v_h,\v_h).
\end{array}
\]
Then by using the last equation, lemmas \ref{propch} and  \ref{lemmapp}, and Assumption \ref{assumpdata}, we deduce that for all $\v_h \in \mathcal{V}_{0h}$, we have
\[
\begin{array}{rcl}
\medskip
(\nabla F(\v_h), \nabla \v_h) &\ge&
a_h(\w_h, \v_h) + c_h(\w_h + {\G}, {\G},\v_h)  - \ell_h(\v_h) + a_h({\G},\v_h) \\
\medskip
&\ge& \ds \frac{\varepsilon_0}{Re} |\v_h |^2_X - c_v (|\v_h|_X + ||{\G} ||_X) \| {\G}\|_X |\v_h |_X\\
&&  -  \varepsilon_0 \| \f \|_{X_0^{'}} | \v_h |_X - \ds (\frac{1}{Re} + \alpha_M C_2^2)|| {\G} ||_X |\v_h |_X.
\end{array}
\]
We get after a simple calculation
\[
(\nabla F(\v_h), \nabla \v_h) \ge  \ds |\v_h |_X \Big ( |\v_h |_X \big( \frac{\varepsilon_0}{Re} -c_v ||{\G} ||_X  \big) -  \ds \varepsilon_0 \| \f \|_{X_0^{'}}  - c_v ||{\G} ||_X^2  - \ds (\frac{1}{Re} + \alpha_M C_2^2)|| {\G} ||_X \Big).
\]
if $\g$ is such that
\begin{equation}\label{rrr6}
|| {\g} ||_{H^{1/2}(\Gamma)^d} < \eta_n = \ds \frac{\varepsilon_0}{2 c_v  Re},
\end{equation}
then we obtain
\[
(\nabla F(\v_h), \nabla \v_h) \ge  \ds |\v_h |_X \Big ( \frac{\varepsilon_0}{2Re} |\v_h |_X  -  \ds \varepsilon_0 \| \f \|_{X_0^{'}}  - c_v  ||{\G} ||_X^2  - \ds  (\frac{1}{Re} + \alpha_M C_2^2)|| {\G} ||_X \Big).
\]
Thus, Brouwer's Fixed-Point Theorem (see for instance \cite{GiraultRaviart}) implies immediately the existence of at least one solution of Problem \eqref{W1h} satisfying the relation
\begin{equation}\label{R1}
|\w_h |_X \le  R_1(\f,\g)= \ds \frac{2 Re}{ \varepsilon_0} ( \varepsilon_0 \| \f \|_{X_0^{'}}  + c_v  || {\g} ||^2_{H^{1/2}(\Gamma)^d}   + \ds  (\frac{1}{Re} + \alpha_M C_2^2)|| {\g} ||_{H^{1/2}(\Gamma)^d}).
\end{equation}
Now we will prove the uniqueness of the solution of Problem \eqref{W1h}. Let $(\w^1_h,p^1_h)$ and $(\w^2_h,p^2_h)$  two solutions of Problem \eqref{W1h} and let $\r_h=\w^1_h-\w^2_h$ and $\xi_h= p^1_h - p^2_h$. Then, $(\r_h, \xi_h)$ is solution of the following problem: $\forall \v_h \in X_{0h}$,
\begin{equation}\label{equatrh}
a_h(\r_h,\v_h) +  c_h(\w^1_h + {\G},\w^1_h + {\G},\v_h) - c_h(\w^2_h + {\G},\w^2_h + {\G},\v_h) - b_h(\v_h, \xi_h) = 0.
\end{equation}
The second and third terms of the last equation can be written as:
\[
\begin{array}{ll}
\medskip
c_h(\w^1_h + {\G},\w^1_h + {\G},\v_h) - c_h(\w^2_h+ {\G},\w^2_h + {\G},\v_h) \\
\medskip
\hspace{2cm}  = c_h(\w^1_h + {\G},\r_h,\v_h)  +  c_h(\w^1_h + {\G},\w^2_h + {\G},\v_h) - c_h(\w^2_h+ {\G},\w^2_h + {\G},\v_h)\\
\medskip
\hspace{2cm} =  c_h(\w^1_h + {\G},\r_h,\v_h)  + d_h(\r_h,\w^2_h + {\G},\v_h) + (\beta_h (|\w^1_h + {\G}| - |\w^2_h + {\G}|) (\w^2_h + {\G}), \v_h).
\end{array}
\]
By taking $\v_h = \r_h$ in \eqref{equatrh}, remarking that $c_h(\w^1_h + {\G},\r_h,\r_h) \ge 0$ and $b_h(\r_h, \xi_h)=0$, using \eqref{rrr6} and \eqref{R1}, using Lemma \ref{lemmapp} for $p=2$ and $4$, and taking into account that $|a|-|b| \le |a-b|$, we get
\[
\begin{array}{rcl}
\medskip
a_h(\r_h,\r_h) &\le& |d_h(\r_h,\w^2_h + {\G},\r_h)| + |(\beta_h (|\w^1_h + {\G}| - |\w^2_h + {\G}|) (\w^2_h + {\G}), \r_h)| \\
\medskip
&\le& (|\w^2_h |_X + ||{\G} ||_X)  \big( c_d + \beta_M C_2 C_4 S_4  \big) |\r_h |_X^2 \\
\medskip
&\le& (c_d + \beta_M C_2 C_4 S_4) ( R_1(\f,\g) + \|\g \|_{H^{1/2}(\Gamma)^d}) |\r_h |_X^2.
\end{array}
\]
Then if $\f$ and $\g$ are sufficiently small such that
\[
\ds R_1(\f,\g) + \|\g \|_{H^{1/2}(\Gamma)^d} < \ds \frac{\varepsilon_0}{Re(c_d + \beta_M C_2 C_4 S_4)},
%\frac{2 c_v Re}{ \varepsilon_0} ( \varepsilon_0 \| \f \|_{X_0^{'}}  + c_{gg} || {\g} ||^2_{H^{1/2}(\Gamma)^d}   + \ds (\frac{1}{Re} + \alpha_M C_2^2)|| {\g} ||_{H^{1/2}(\Gamma)^d})  < \ds \frac{\varepsilon_0}{2Re},
\]
which gives
\[
 \varepsilon_0 \| \f \|_{X_0^{'}}  + c_v || {\g} ||^2_{H^{1/2}(\Gamma)^d}   + \ds (\frac{\varepsilon_0 + 2}{2Re }  + \alpha_M C_2^2)|| {\g} ||_{H^{1/2}(\Gamma)^d}  < \ds \frac{\varepsilon_0^2}{2Re^2(c_d + \beta_M C_2 C_4 S_4)},
\]
then we get
\[
a_h(\r_h,\r_h) < \ds \frac{\varepsilon_0}{Re} |\r_h |^2_X.
\]
We deduce that $\r_h=0$ and then $\r^1_h = \r^2_h$. Furthermore, Relation \eqref{equatrh} deduces that $b_h(\v_h,\xi_h)=0$ which gives with the inf-sup condition \eqref{infsuph1h} that $\xi_h=0$. Thus we get the uniqueness of the solution.\hfill$\Box$\\

In \cite{BDFtheornum1} the authors show an {\it a priori} error estimate corresponding to their discrete scheme by using the Brezzi-Rappaz Theorem. Even the discrete formulation \eqref{W1h} contain a supplementary term of stabilisation, the corresponding {\it a priori} error estimate can similarly be established and we get the following error between the exact solution $(\u,p)$ and the discrete one $(\u_h, p_h)$: if the exact solution $(\w,p)$ of Problem \eqref{W1} is in $H^2(\Omega)^d \times H^1(\Omega)$, we have,
\[
\begin{array}{rcl}
\medskip
\|\u -\u_h \|_X  + \|p - p_h \|_{L^2(\Omega)} &\le& C h (\| \w \|_{H^2(\Omega)^d} + \|p \|_{H^1(\Omega)}) + \| \G - \mathcal{I}_{X_h} \G \|_X \\
&& + C h \max (\|\varepsilon \|_{W^{1,\infty}(\Omega)}, \| \varepsilon \|_{W^{2,r}(\Omega)}).
\end{array}
\]
{\rmq The exact and numerical solutions of Problems \eqref{W1} and \eqref{W1h} must be sufficiently small to get the corresponding uniqueness (see Theorems \ref{exuniqtheo} and  \ref{exuniqnum}). Then we can establish the {\it a priori} error estimate by using the classical method which considers the difference between the numerical Problem \eqref{W1h} and the exact one \eqref{W1} (for $\v=\v_h$) and which uses adapted interpolation operators: Operator $P_h$ (see \cite{GiraultLions}, page 35) for the velocity and Scott-Zhang operator $\mathcal{F}_h$ for the pressure (see \cite{CiarletP})).
}
\subsection{Iterative algorithms}
In order to approximate the solution of the non-linear problem \eqref{W1h}, we introduce the following iterative algorithm. \\

For a given initial guess $\w_h^0 \in X_{0h}$ and having $\w_h^i$ at each iteration $i$, we compute $(\w_h^{i+1}, p^{i+1}_h)$ solution~of:
\begin{equation}\label{Whi}
\left\{
\begin{array}{ll}
\medskip
\forall \v_h \in X_{0h}, \quad \ds   a_h(\w_h^{i+1},\v_h) + c_h(\w_h^i + {\G},\w_h^{i+1},\v_h)  - b_h(\v_h,p_h^{i+1}) \\
\medskip
\hspace{6cm} = \ell_h(\v_h)  -a_h({\G},\v_h) - c_h(\w_h^i + {\G}, {\G},\v_h)  ,\\
\forall q_h\in M_{0h}, \quad \ds b_h(\w_h^{i+1}, q_h) = 0.
\end{array}
\right.
\end{equation}
In the following, we investigate the existence and uniqueness, and the convergence of the solution of Scheme (\ref{Whi}). We begin by proving the existence and uniqueness.

{\thm \label{existence}
Problem  \eqref{Whi} admits a unique solution in $X_{0h} \times M_{0h}$.
}

\noindent \textbf {Proof.}
To prove the existence and  uniqueness of the solution of Problem $(\ref{Whi})$ which is a square linear system in finite dimension, it suffices to show the uniqueness. For a given $\w_h^i$, let $(\w_{h1}^{i+1},p_{h1}^{i+1})$ and $(\w_{h2}^{i+1},p_{h2}^{i+1})$ two different solutions of Problem $(\ref{Whi})$ and let $\w_h= \w_{h1}^{i+1} - \w_{h2}^{i+1}$ and $\xi_h = p_{h1}^{i+1} - p_{h2}^{i+1}$, then $(\w_h,\xi_h)$ is the solution of the following problem:
\begin{equation}\label{rrr3}
\forall \v_h \in \mathcal{V}_{0h}, \quad \ds  a_h(\w_h,\v_h) + c_h(\w_h^i + {\G},\w_h,\v_h) =  0.
\end{equation}
By taking $\v_h = \w_h$, remarking that $c_h(\w_h^i + {\G},\v_h,\v_h)\ge 0$ (Lemma \ref{propch}), we get
\[
 \alpha_m \|\w_h \|_{L^2(\Omega)^d}^2 + \ds \frac{\varepsilon_0}{Re} |\w_h |_X^2 \le 0,
\]
and then $\w_h = 0$. The inf-sup condition \eqref{infsuph1h} gives $\xi_h=0$. Thus, we get the existence and the uniqueness of the solution of Problem \eqref{Whi}.
$\hfill\Box$\\

Next, the following theorem shows that the solution of Problem \eqref{Whi} is bounded for small values of the data.
{\thm \label{boundu1} We consider Problem \eqref{Whi}. For a given initial guess $\w_h^0$ such that
\[
|\w_h^0 |_X \le F_1(\f,\g) = \ds \frac{2 Re}{\varepsilon_0} \big( \|\f \|_{X_0^{'}} +  (\frac{1}{Re} + \alpha_M S_2) \|\g \|_{H^{1/2}(\Gamma)^d} + c_v \|\g \|_{H^{1/2}(\Gamma)^d}^2   \big)
\]
and if the data $\g$ is small such that
\[
\|\g \|_{H^{1/2}(\Gamma)^d} \le \ds \frac{\varepsilon_0}{2c_v Re},
\]
then the solution of Problem \eqref{Whi} is bounded as following
\begin{equation}\label{equatt2}
|\w_h^{i+1} |_X \le F_1(\f,\g).
\end{equation}
}
\noindent \textbf {Proof.}
To show Relation \eqref{equatt2},  we consider the first equation of (\ref{Whi}) with $\v_h = \w^{i+1}_h$ and we get:
\begin{equation}
\ds   a_h(\w_h^{i+1},\w_h^{i+1}) + c_h(\w_h^i + {\G},\w_h^{i+1},\w_h^{i+1})
= \ell_h(\v_h) -a_h({\G},\w_h^{i+1}) - c_h(\w_h^i + {\G}, {\G} ,\w_h^{i+1}) .
\end{equation}
By remarking that $c_h(\w_h^i + {\G},\w_h^{i+1},\w_h^{i+1}) \ge 0 $ and by using Lemma \ref{propch} we have
\begin{equation}
\begin{array}{rcl}
\medskip
\ds
\frac{\varepsilon_0}{Re} | \w_h^{i+1} |_X^2
&\le& \Big| \ell_h(\w_h^{i+1}) - a_h({\G},\w_h^{i+1}) -c_h(\w_h^i + {\G},{\G},\w_h^{i+1}) \Big|\\
\medskip
&\le& \ds |\w_h^{i+1} |_X  \big( \|\f \|_{X_0^{'}} +  (\frac{1}{Re} + \alpha_M S_2) \|{\G} \|_X \big) + c_v  \|{\G} \|_X (| \w_h^i |_X + \|{\G} \|_X)  | \w_h^{i+1} |_X.
\end{array}
\end{equation}
Then we get
\[
\frac{\varepsilon_0}{Re} | \w_h^{i+1} |_X \le \ds \|\f \|_{X_0^{'}} +  (\frac{1}{Re} + \alpha_M S_2) \|{\G} \|_X + c_v \|{\G} \|_X^2  + c_v   | \w_h^i |_X \|{\G} \|_X.
\]
We denote by
\[
F_1(\f,\g) = \ds \frac{2 Re}{\varepsilon_0} \big( \|\f \|_{X_0^{'}} +  (\frac{1}{Re} + \alpha_M S_2) \|\g \|_{H^{1/2}(\Gamma)^d} + c_v \|\g \|_{H^{1/2}(\Gamma)^d}^2   \big),
\]
then we have
\begin{equation}\label{inducti}
\frac{\varepsilon_0}{Re} | \w_h^{i+1} |_X \le \ds \frac{\varepsilon_0}{2Re} F_1(\f,\g)  + c_v   | \w_h^i |_X \|\g \|_{H^{1/2}(\Gamma)^d}.
\end{equation}
We show Bound \eqref{equatt2} by induction on $i$. We suppose that the initial guess $\w_h^0$ is such that
\[
|\w_h^0 |_X \le \ds F_1(\f,\g).
\]
and we suppose that
\[
|\w_h^i |_X \le \ds F_1(\f,\g).
\]
We are in one of the following two situations:
\begin{itemize}
\item  We have $| \w_h^{i+1} |_X \le | \w_h^i |_X$, we deduce immediately that
\[
| \w_h^{i+1} |_X \le F_1(\f,\g).
\]
\item  We have $| \w_h^{i+1} |_X > | \w_h^i |_X$, Equation \eqref{inducti} gives
\[
\frac{\varepsilon_0}{Re} | \w_h^{i+1} |_X \le \ds \frac{\varepsilon_0}{2Re} F_1(\f,\g)  + c_v   | \w_h^{i+1} |_X \|\g \|_{H^{1/2}(\Gamma)^d}.
\]
As the data $g$ is small such that
\[
\|\g \|_{H^{1/2}(\Gamma)^d} \le \ds \frac{\varepsilon_0}{2c_v Re},
\]
we get
\[
| \w_h^{i+1} |_X \le F_1(\f,\g).
\]
\end{itemize}
We deduce finally the bound \eqref{equatt2}. \hfill$\Box$

The next theorem shows the convergence of the solution $(\u_h^i, p_h^i)$ of Problem $(\ref{Whi})$ in $X \times L^2(\Omega)$ to the solution
$(\u_h, p_h)$ of Problem $(\ref{W1h})$.
{\thm \label{converg1} Under the assumptions of Theorems \ref{exuniqnum} and \ref{boundu1}, and if the data $\f$ and $\g$ are sufficiently small (see \eqref{sufsmall}), then the sequence of solutions $(\u_h^{i+1}, p_h^{i+1})$ of Problem $(\ref{Whi})$ converges in $L^2(\Omega)^d \times L^2(\Omega)$ to the solution
$(\u_h, p_h)$ of Problem $(\ref{W1h})$. }

\noindent \textbf {Proof.} We consider the difference between the
equations (\ref{Whi}) and (\ref{W1h}) and we obtain for $\v_h^{i+1} = \v_h =  \w_h^{i+1} - \w_h$,
\begin{equation}\label{relat22}
\ds  a_h(\v_h^{i+1},\v_h^{i+1}) + c_h(\w_h^i+{\G},\w_h^{i+1}+{\G} ,\v_h^{i+1}) - c_h(\w_h+{\G},\w_h+{\G},\v_h^{i+1})  = 0.
\end{equation}
We insert $\w_h$ in the second term of the last equation and we get
\[
a_h(\v_h^{i+1},\v_h^{i+1}) + c_h(\w_h^i+{\G},\v_h^{i+1},\v_h^{i+1}) + c_h(\w_h^i+{\G},\w_h+{\G},\v_h^{i+1}) - c_h(\w_h+{\G},\w_h+{\G},\v_h^{i+1})  = 0.
\]
By using the properties of the form $c_h$, the third and fourth terms of the last equation can be written as:
\[
\begin{array}{ll}
\medskip
c_h(\w_h^i+{\G},\w_h+{\G},\v_h^{i+1}) - c_h(\w_h+{\G},\w_h+{\G},\v_h^{i+1})
\\
\medskip
\hspace{2cm}
=  (\varepsilon_h (|\w_h^i+{\G}| - |\w_h+{\G} |), \w_h+{\G},\v_h^{i+1})
%\\
%\medskip
%\hspace{5cm}
+ d_h(\v_h^i, \w_h+{\G}, \v_h^{i+1}).
\end{array}
\]
Then we get
\[
\begin{array}{ll}
\medskip
\ds  a_h(\v_h^{i+1},\v_h^{i+1}) + c_h(\w_h^i+{\G},\v_h^{i+1},\v_h^{i+1}) = -(\varepsilon_h (|\w_h^i+{\G}| - |\w_h+{\G} |), \w_h+{\G},\v_h^{i+1}) \\
\medskip
\hspace{7cm} - d_h(\v_h^i, \w_h+{\G}, \v_h^{i+1}).
\end{array}
\]
By remarking that $c_h(\w_h^i+{\G},\v_h^{i+1},\v_h^{i+1})\ge 0$, We obtain the following bound:
\[
\begin{array}{rcl}
\medskip
\ds  \alpha_m \|\v_h^{i+1} \|_{L^2(\Omega)^d}^2 +  \frac{\varepsilon_0}{Re} |\v_h^{i+1} |_X^2 &\le& S_2 S_4 C_4 (|\w_h^i - \w_h|_X ) (\| \w_h \|_X + \|{\G} \|_X) |\v_h^{i+1} |_X \\
\medskip
&& + c_d (\|\w_h \|_X + \| {\G} \|_X) |\v_h^i |_X |\v_h^{i+1} |_X\\
&\le& c_1 (| \w_h |_X + \|{\G} \|_X) |\v_h^i |_X |\v_h^{i+1} |_X.
\end{array}
\]
Relation \eqref{majwh} gives
\begin{equation}\label{relatconv}
\begin{array}{ll}
\ds  \alpha_m \|\v_h^{i+1} \|_{L^2(\Omega)^d}^2 +\frac{\varepsilon_0}{2Re} |\v_h^{i+1} |_X \le c_1 (R_1(\f,\g) + \|\g \|_{H^{1/2}(\Gamma)^d})  |\v_h^i |_X |\v_h^{i+1} |_X.
\end{array}
\end{equation}
If the data $\f$ and $\g$ are sufficiently small such that
\begin{equation}\label{sufsmall}
c_1 (R_1(\f,\g) + \|\g \|_{H^{1/2}(\Gamma)^d}) < \ds \frac{\varepsilon_0}{2Re},
\end{equation}
then we get
\[
|\v_h^{i+1} |_X \le |\v_h^i |_X,
\]
and we deduce the convergence of the sequence $(\w_h^{i+1} - \w_h)$ in $X$ and then the convergence of the sequence $\w_h^{i}$ in $X_0$. By taking the limit of (\ref{relatconv}), we get
\[
\underset{{\small i \rightarrow +\infty} }{\lim} \Big( {  || \w_h^{i+1} - \w_h ||^2_{L^2(\Omega)}}  \Big) \le 0.
\]
We deduce then that ${  || \w_h^{i+1} - \w_h ||_{L^2(\Omega)} }$ converges to $0$ and $\w_h^{i+1}$ converges to $\w_h$ in $L^2(\Omega)^d$.\\

For the convergence of the pressure, We take the difference between the equations (\ref{Whi}) and (\ref{W1h}) and we obtain for all $\v_h \in X_{0h}$ the equation
\[
\begin{array}{rcl}
\medskip
\ds b_h(\v_h, p_h^{i+1} - p_h) &=& \ds  a_h(\w_h^{i+1} - \w_h, \v_h) +
c_h(\w_h^i + {\G}, \w_h^{i+1} + {\G},\v_h) - c_h(\w_h + {\G}, \w_h + {\G},\v_h)\\
&=& a_h(\w_h^{i+1} - \w_h, \v_h) + c_h(\w_h^i + {\G}, \w_h^{i+1} -\w_h,\v_h) \\
\medskip
&& + c_h(\w_h^i + {\G}, \w_h + {\G},\v_h) - c_h(\w_h + {\G}, \w_h + {\G},\v_h)\\
&=& a_h(\w_h^{i+1} - \w_h, \v_h) + c_h(\w_h^i + {\G}, \w_h^{i+1} -\w_h,\v_h) +d_h(\w_h^i - \w_h, \w_h + {\G},\v_h) \\
\medskip
&& + (\beta_h (|\w_h^i + {\G}| - |\w_h + {\G}|) (\w_h + {\G}), \v_h).
\end{array}
\]
We get then the following bound:
\[
\begin{array}{ll}
\medskip
\ds \frac{\Big| b_h(\v_h, p_h^{i+1} - p_h) \Big| }{||\v_h ||_X} \le  \ds \frac{1}{Re} |\w_h^{i+1} - \w_h|_X + \alpha_M \|\w_h^{i+1} - \w_h \|_{L^2(\Omega)^d} + c_v \|\w_h^i +{\G} \|_X |\w_h^{i+1} - \w_h|_X \\
\medskip
\hspace{4cm} + c_v \|\w_h +{\G}\|_X |\w_h^i - \w_h|_X + \beta_M S_2S_4 C_4 \|\w_h + {\G} \|_X |\w_h^i - \w_h|_X.
\end{array}
\]
Owning the inf-sup condition \eqref{infsuph1h} and taking the limit of the last equation deduces the strong convergence of the subsequence $p_h^i$ to $p_h$ in $L^2(\Omega)$.\hfill$\Box$\\
\section{A posteriori error estimate}
\noindent We start this section by introducing some additional notations and properties that will be useful in order to establish an {\it a posteriori} estimate.\\

%---------------------redefine notations
\noindent For any element $\kappa$ in $\mathcal{T}_{h}$, we denote by:\\
$\bullet\;$  $\varepsilon_{\kappa}$ the set of edges (when $ d=2$) or faces (when $ d=3$) of $\kappa$ that are not contained in  $\Gamma$,\\
$\bullet\;$
$h_{\kappa}$ the diameter of the element $\kappa$ and $h_{e}$ the  diameter of edge (or face) $e$,\\
$\bullet\;$ $[\cdot]_{e}$ the jump through $e$ on each edge (or face) $e$ on  $\varepsilon_{\kappa},$ \\
$\bullet\;$ ${ \n}_{\kappa}$  stands for the unit outward normal vector to $\kappa $ on $\partial \kappa$.\\

\noindent We denote also by $\mathcal{E}_h$ the set of all the edges that are not containing in $\Gamma$. In other term, $\mathcal{E}_h$ is the set of the interior edges of the mesh.\\

%-------------------
\noindent We introduce the  following inverse  inequalities, (see \cite{ErnGuerm}, page 75):
For any number $p\geq 2$, for any dimension $d$, and for any non negative integer $r$, there exist  constants $c_I^0(p)$ such that for any polynomial function $v_h$ of degree $r$  on an element $\kappa$ of $\mathcal{T}_{h}$,
\begin{equation} \label{dk}
 \|v_h\|_{L^p(\kappa)}\leq c_I^0 (p) h_{\kappa}^{\frac{d}{p}-\frac{d}{2}}\|v_h\|_{L^2(\kappa)}.
\end{equation}

\noindent We now recall the following definitions and properties (see R. Verf\"urth, \cite[Chapter 1]{Verfurth2013}): for an element $\kappa$ of $\mathcal{T}_h$, the bubble function $\psi_\kappa$ (resp. $\psi_e$ for the face $e$) is defined as  the product of the $d+1$ barycentric coordinates associated with the vertices of $\kappa$ (resp. of the $d$ barycentric coordinates associated with the vertices of~$e$). We also consider a lifting operator ${\mathcal{L}}_{e}$ defined on polynomials on $e$ vanishing on $\partial e$ into polynomials on the at most two elements $\kappa$ containing $e$ and vanishing on $\partial \kappa \setminus e $, which is constructed by affine transformation from a fixed operator on the reference element.
\begin{propri}\label{psi}
 Denoting by $P_r(\kappa)$ the space of polynomials of degree smaller than or equal to $r$ on $\kappa$, the following properties hold:
\begin{equation}\label{eq:ineg_verfurth}
\forall v \in P_r(\kappa), \qquad
\begin{cases}
c ||v||_{L^2(\kappa)} \le ||v \psi^{1/2}_{\kappa} ||_{L^2(\kappa)}
\le c' ||v||_{L^2(\kappa)}, &\\
|v|_{H^1(\kappa)} \le c h_{\kappa}^{-1} ||v ||_{L^2(\kappa)}.&
\end{cases}
\end{equation}
\end{propri}
\begin{propri} \label{psii} Denoting by $P_r(e)$ the space of polynomials of degree smaller than or equal to $r$ on $e$, we have
$$
\forall\; v  \in P_r(e),\qquad
c\Vert v \Vert_{L^2(e)}\leq \Vert v\psi_{e}^{1/2}
\Vert_{L^2(e)}\leq c'\Vert  v \Vert_{L^2(e)},
$$
and, for all polynomials $v$ in $P_r(e)$ vanishing on $\partial e$, if $\kappa$ is an element which contains $e$,
$$ \Vert {\mathcal{L}}_{e}v \Vert_{L^2(\kappa)}+h_{e}\mid
{\mathcal{L}}_{e}v \mid_{H^1(\kappa)}\leq ch^{1/2}_{e}\Vert  v
\Vert_{L^2(e)}.$$
\end{propri}
We also introduce the Cl\'ement type regularization operator $C_h$  which has the following properties, see \cite [section IX.3]{BMR} and \cite{clement}: For any function $\w$ in $H^1(\Omega)^d$, $C_h \w$  belongs to the continuous affine finite element space and satisfies for any $\kappa$ in $\mathcal{T}_h$ and $e$ in $\varepsilon_{\kappa}$,
\begin{equation}\label{clement1}
\begin{array}{rcl}
||\w-C_h \w||_{L^2(\kappa)^d} &\le& ch_\kappa||\w||_{H^1(\Delta_\kappa)^d} \quad
\mbox{ and } \quad
||\w-C_h \w||_{L^2(e)^d} \le ch^{1/2}_e ||\w||_{H^1(\Delta_e)^d},
\end{array}
\end{equation}
where $\Delta_\kappa$ and $\Delta_e$ are the following sets:
$$
\Delta_\kappa = \bigcup \Big\{ \kappa' \in \mathcal{T}_h;\kappa'\cap \kappa \ne 0 \Big\} \quad \mbox{ and  } \quad \Delta_e = \bigcup \Big\{ \kappa' \in \mathcal{T}_h;\kappa'\cap e \ne 0 \Big\}. $$
Note that we use the variant of $C_h$ which ensures that $C_h \w$ belongs to $H^1_0(\Omega)^d$ (see \cite{clement}). Furthermore more we have the continuity relation: for all $\v\in X_0$, $C_h \v \in X_{0h}$ and
\begin{equation}\label{Chcont}
 |C_h \v |_{X} \le c |\v|_X.
\end{equation}
%

%----------------construction of indicators
In this section, we specify the operator $\mathcal{I}_{X_h}$ to be the Cl\'ement operator $C_h$. We also introduce averaged values for data; this approximations will be useful to prove the optimally of the indicators. Let $\f_h, \alpha_{\varepsilon_h}$ and $\beta_{\varepsilon_h}$ be the piecewise constant approximations of the data $\f,\alpha$ and $\beta$:
$$
\f_h |\kappa= \ds \frac{1}{|\kappa|} \int_{\kappa}  \f (\x) d\x, \quad \alpha_{\varepsilon_h} |\kappa= \ds \frac{1}{|\kappa|} \int_{\kappa} \alpha (\varepsilon(\x)) d\x, \quad \beta_{\varepsilon_h} |\kappa= \ds \frac{1}{|\kappa|} \int_{\kappa} \beta (\varepsilon(\x)) d\x.
$$
Next, we distinguish the discretization and linearization estimators. For this, we first write the residual equation. The difference between  \eqref{W1} and  \eqref{Whi} gives the following relations for all $\v \in X_0$ and all $\v_h \in X_{0h}$:
\begin{equation}\label{equation1}
\begin{array}{ll}
\medskip
\Big( a(\w,\v) + c(\w+{\G},\w+{\G},\v)  - b(\v,p) \Big) - \Big(a(\w_h^{i+1},\v) + c(\w^{i+1}_h+{\G},\w_h^{i+1}+{\G},\v) - b(\v,p_h^{i+1})\Big) \\
\medskip
\hspace{1cm}= \mathcal{F}_\varepsilon(\v) +\mathcal{F}_{it}(\v) + \mathcal{F}(\v_h - \v),
\end{array}
\end{equation}
where
\begin{equation}\label{FFepsi}
\begin{array}{rcl}
\medskip
\mathcal{F}_\varepsilon(\v) &=& \big( (\ell-\ell_h)(\v) -(a-a_h)({\G},\v) \big) \\
 \medskip
&& + \big( (a_h-a)(\w_h^{i+1},\v) + (c_h-c)(\w_h^{i+1}+{\G},\w_h^{i+1}+{\G},\v) - (b_h-b)(\v,p_h^{i+1}) \big),\vspace{.5cm}
\end{array}
\end{equation}
\begin{equation}\label{Fit}
\hspace{-6cm} \mathcal{F}_{it} (\v) = -c_h(\w^{i+1}_h+{\G},\w_h^{i+1}+{\G},\v) + c_h(\w^i_h+{\G},\w_h^{i+1}+{\G},\v), \vspace{.5cm}
\end{equation}
\begin{equation}\label{FF}
\begin{array}{rcl}
\medskip
\mathcal{F}(\v_h - \v) &=& \big( (\ell_h(\v-\v_h) -a_h({\G},\v-\v_h) \big)\\
\medskip
&& + \big( a_h(\w_h^{i+1},\v_h-\v) + c_h(\w_h^i+{\G},\w_h^{i+1}+{\G},\v_h-\v)
- b_h(\v_h-\v,p_h^{i+1}) \big),
\end{array}
\end{equation}
and for all $q\in M$:
\begin{equation}\label{equation2}
 b(\w- \w^{i+1}_h,q)= - (b-b_h)(\w_h^{i+1},q) - b_h(\w_h^{i+1},q).
\end{equation}
By Adding and subtracting $\f_h, \alpha_{\varepsilon_h}, \beta_{\varepsilon_h}$ and ${\G}_h = C_h \G$ in Equation \eqref{FF}, and using the the Green formula on each $\kappa \in \mathcal{T}_h$, we obtain the following formula:
\begin{equation}\label{FFepsivvh}
\begin{array}{rcl}
\medskip
\mathcal{F} (\v_h - \v) &=& \ds ( \varepsilon_h (\f - \f_h) , \v - \v_h) + c_h(\w_h^i+{\G},{\G}-{\G}_h,\v_h-\v) +
d_h({\G} - {\G}_h, \w_h^{i+1} + {\G}_h, \v_h - \v)\\
\medskip
&& + (\beta_h (|\w_h^i + {\G}| - |\w_h^i + {\G}_h|)( \w_h^{i+1} + {\G}_h), \v_h - \v) -a_h({\G}-{\G}_h,\v-\v_h) \\
\medskip
&& + \ds \sum_{\kappa \in \mathcal{T}_h} \Big\{ \int_\kappa \Big(\varepsilon_h \f_h + \frac{1}{Re}  \div (\varepsilon_h  \nabla (\w_h^{i+1} + {\G}_h)  ) - \alpha_{\varepsilon_h} (\w_h^{i+1}+{\G}_h) \\
\medskip
&& \qquad \qquad \ds - \varepsilon_h (({\G}_h + \w_h^i)\cdot \nabla ) ({\G}_h+\w_h^{i+1}) -
\frac{1}{2} \div(\varepsilon_h ({\G}_h + \w_h^i)) ({\G}_h+ \w_h^{i+1}) \\
\medskip
&&
 \qquad \qquad  - \beta_{\varepsilon_h} |{\G}_h + \w_h^i|  ({\G}_h+ \w_h^{i+1}) - \varepsilon_h \nabla p_h^{i+1}\Big) (\v - \v_h) d\x\\
\medskip
&& \qquad \quad - \ds \frac{1}{2} \sum_{e \in \varepsilon_{\kappa} \cap \mathcal{E}_h} \int_{e}  \ds [(\frac{1}{Re} \varepsilon_h \nabla (\w_h^{i+1} + {\G}_h) - p_h^{i+1} \mathbb I)(\sigma)\cdot {\n}] (\v - \v_h) d \sigma
\Big\}.
\end{array}
\end{equation}
The term $\mathcal{F}_\varepsilon$ will be bounded by data errors as the numerical solutions are bounded in $X_0 \times M$. The term $\mathcal{F}_{it}$, $\mathcal{F}$ and Equation \eqref{equation2} allow us to define the following  local linearization indicator $(\eta_{i,\kappa}^L) $ and local discretization indicator $(\eta_{i,\kappa}^{D})$:

\begin{equation}\label{indiclin}
(\eta_{i,\kappa}^L)=   \|\w_h^{i+1}- \w_h^{i} \|_{H^1(\kappa)^d},\vspace{.02cm}
\end{equation}
\begin{equation}\label{inddis}
\begin{array}{rcl}
(\eta_{i,\kappa}^{D}) &=& \ds h_{\kappa } || \varepsilon_h \f_h + \frac{1}{Re}  \div (\varepsilon_h  \nabla (\w_h^{i+1} + {\G}_h)  ) - \alpha_{\varepsilon_h} (\w_h^{i+1}+{\G}_h) \\
\medskip
&& \qquad \qquad \ds - \varepsilon_h (({\G}_h + \w_h^i)\cdot \nabla ) ({\G}_h + \w_h^{i+1}) -
\frac{1}{2} \div(\varepsilon_h ({\G}_h + \w_h^i)) ({\G}_h + \w_h^{i+1}) \\
\medskip%
&& \qquad \qquad - \beta_{\varepsilon_h} (|{\G}_h|+ |\w_h^i|) ({\G}_h + \w_h^{i+1})  - \varepsilon_h \nabla p_h^{i+1} ||_{L^2(\kappa)^d}\\
\medskip
&& + \ds\frac{1}{2} \sum_{e \in \varepsilon_{\kappa}\cap \mathcal{E}_h} h^{1/2}_{e} ||[(\frac{1}{Re} \varepsilon_h \nabla (\w_h^{i+1} + {\G}_h) - p_h^{i+1} \mathbb I)(\sigma)\cdot {\n}] ||_{L^2(e)^d}
 + ||\div (\varepsilon_h \w^{i+1}_h)||_{L^2(\kappa)}.
\end{array}
\end{equation}
We now are in position to bound the error between the exact and numerical solutions with the indicators up to the errors on the data.
{\thm \label{upperb}
Let  $\W=(\w,p) $ be a non singular solution of Problem \eqref{W1}. Then, there exists a neighborhood $\Theta $ of $ \W$ in  $X_0 \times M$ such that any solution  $\W^{i+1}_h=(  \w^{i+1}_h,  \p^{i+1}_h)  \in X_{0h} \times M_{0h} $ in  $ \Theta$ of Problem \eqref{Whi} satisfies the following {\it a posteriori} error estimate:
 \[
\begin{array}{l}
\medskip
||\w- \w^{i+1}_h||_{X} +
||p  -  \p^{i+1}_h||_{M} \le \ds C \Big( \sum_{\kappa \in \mathcal{T}_h}  \big( h^2_\kappa  || \f -  \f_h ||^2_{L^2(\kappa)^d} + | {\G} - {\G}_h |^2_{H^1(\kappa)^d}  + ||\varepsilon - \varepsilon_h ||^2_{L^\infty(\kappa)}\\
\hspace{.5cm} \ds+ ||\nabla \varepsilon - \nabla \varepsilon_h ||^2_{L^3(\kappa)} + \| \alpha (\varepsilon) - \alpha_{\varepsilon_h}\|^2_{L^2(\kappa)} + \| \beta (\varepsilon) - \beta_{\varepsilon_h}\|^2_{L^6(\kappa)} \big) + \sum_{\kappa \in \mathcal{T}_h} \bigl(\eta_{i,\kappa}^D)^2  + \sum_{\kappa \in \mathcal{T}_h}  \bigl(\eta_{i,\kappa}^L)^2  \Big)^{1/2},
\end{array}
\]
where $C$ is a constant depending on $\W$ but independent of $h$.
}

\noindent \textbf {Proof.}
\noindent  Let $\W=(\w,p) $ be a non singular solution of Problem \eqref{W1} and  $ \W^{i+1}_h=(  \w^{i+1}_h,  \p^{i+1}_h)  \in X_{0h} \times M_{0h} $ be the solution of the iterative problem  \eqref{Whi}. Having proved Proposition \ref{isom}, and owing to  \cite{Poussin} and \cite [Prop. 2.2]{Verfurth}, there exists a neighborhood $\Theta$ of $ \W$ in  $X_0 \times M$ such that if $\W^{i+1}_h$  is in $\Theta$, we have the following bound:
% all $\v_h \in X_{0h}$,
\begin{equation*}
\begin{array}{rcl}
\medskip
|| \W- \W^{i+1}_h ||_{X_0 \times M} & \le &  c || F(\W) - F( \W^{i+1}_h)||_{X_0 \times M}\\
& \le &   c ||  \W^i_h+ \mathcal S (\mathcal{G}(  \w^{i+1}_h),0)- \W- \mathcal S (\mathcal{G}(\w),0)||_{X_0 \times M}\\
& \le &  c || \mathcal S \Big(  \mathcal S^{-1}(  \W^{i+1}_h)  + (\mathcal{G}(  \w^{i+1}_h), 0)-\mathcal S^{-1}(\W) - (\mathcal{G}( \w), 0) \Big)||_{X_0 \times M}.
\end{array}
\end{equation*}
\noindent Using \eqref{stokes3}, we have:
\begin{equation*}
\begin{array}{rcl}
\medskip
|| \W- \W^{i+1}_h ||_{X \times M}& \le &    c' ||  \mathcal S^{-1}(  \W^{i+1}_h)  + (G(  \w^{i+1}_h),0) -\mathcal S^{-1}( \W)  - (G(\w), 0) ||_{X_0^{'} \times M'}\\
\end{array}
\end{equation*}
\noindent Relations \eqref{equation1} and \eqref{equation2} allow us to get: for all $\v_h \in X_h$,
\begin{equation} \label{majorationerreur}
\begin{array}{rcl}
\medskip \medskip
 || \W- \W^{i+1}_h ||_{X_0 \times M}& \le &  \ds  \bar{C}_1 \,  \Big(  \underset { \v \ne \0 } {\underset { \v \in X_0} {\sup} } \frac{{\mathcal F}_\varepsilon (\v)}{|\v |_{X_0}} + \underset { \v \ne \0 } {\underset { \v \in X_0} {\sup} }  \frac{{\mathcal F_{it}}(\v)}{ |(\v,q) |_{X_0 \times M}} + \underset { \v \ne \0 } {\underset { \v \in X_0} {\sup} }  \frac{{\mathcal F}(\v_h - \v)}{ |(\v,q) |_{X_0 \times M}}   \\
 && \ds \qquad \quad  +\underset { q \ne 0 } {\underset { q \in M} {\sup} } \ds \frac{   |(b-b_h)(\w_h^{i+1},q) + b_h(\w_h^{i+1},q)|}   {||q||_{M}}
 \Big).
 \end{array}
 \end{equation}
 We will bound every term of the right hand side of \eqref{majorationerreur}.\\
 \noindent We begin with the first one where $\mathcal{F}_\varepsilon (\v)$ is given by \eqref{FFepsi}. The first term of the right hand side of \eqref{FFepsi} gives by using the discrete Cauchy-Schwartz inequality:
 \begin{equation}\label{term1}
 \begin{array}{rcl}
 \medskip
 \big| (\ell-\ell_h)(\v) \big| &=& \ds \big| \sum_{\kappa \in \mathcal{T}_h} \int_\kappa (\varepsilon - \varepsilon_h) \f \v d\x \big| \\
 \medskip
 &\le& \ds \sum_{\kappa \in \mathcal{T}_h} \| \varepsilon - \varepsilon_h \|_{L^\infty(\kappa)} \| \f \|_{L^2(\kappa)^d} \|\v \|_{L^2(\kappa)^d}\\
 &\le& c_1 \| \f \|_{L^2(\Omega)^d} \Big(  \ds \sum_{\kappa \in \mathcal{T}_h} \| \varepsilon - \varepsilon_h \|^2_{L^\infty(\kappa)}  \Big)^{1/2} |\v |_X.
 \end{array}
 \end{equation}
The second term can be treated as follow:
\begin{equation}\label{term2}
\begin{array}{rcl}
\medskip
\big| (a-a_h)({\G}, \v) \big| &\le& \ds \sum_{\kappa \in \mathcal{T}_h} \int_\kappa \frac{1}{Re} \|\varepsilon - \varepsilon_h\|_{L^\infty(\kappa)}  |{\G}|_{H^1(\kappa)^d} |\v|_{H^1(\kappa)^d} d\x  \\
\medskip
&& \qquad \qquad + \| \alpha (\varepsilon) - \alpha_{\varepsilon_h}\|_{L^2(\Omega)} \|{\G}\|_{L^4(\Omega)^d} \|\v\|_{L^4(\Omega)^d}  \\
\medskip
&\le& \ds c_2 \| {\G}\|_{X} \Big[  \Big( \sum_{\kappa \in \mathcal{T}_h} \| \varepsilon - \varepsilon_h \|^2_{L^\infty(\kappa)} + \sum_{\kappa \in \mathcal{T}_h} \| \alpha (\varepsilon) - \alpha_{\varepsilon_h}\|^2_{L^2(\kappa)} \Big)^{1/2} |\v |_X.
\end{array}
\end{equation}
The third term can be treated exactly like the second one and we get
\[
\begin{array}{rcl}
\medskip
\big| (a-a_h)(\w_h^{i+1}, \v) \big|
&\le& \ds c_3 \| \w_h^{i+1} \|_{X} \Big[  \Big( \sum_{\kappa \in \mathcal{T}_h} \| \varepsilon - \varepsilon_h \|^2_{L^\infty(\kappa)} + \sum_{\kappa \in \mathcal{T}_h} \| \alpha (\varepsilon) - \alpha_{\varepsilon_h}\|^2_{L^2(\kappa)} \Big)^{1/2} |\v |_X.
\end{array}
\]
For the second term we have by using the fact that $\w_h^{i+1}$ is bounded in $X$ and the Cauchy-Schwartz inequality:
\[
\begin{array}{ll}
\medskip
\big|  (c-c_h)(\w_h^{i+1}+{\G},\w_h^{i+1}+{\G},\v) \big| \le \big|  (d-d_h)(\w_h^{i+1}+{\G},\w_h^{i+1}+{\G},\v) \big|  \\
\medskip
\hspace{5cm} + \big|  ((\beta(\varepsilon_h) - \beta_{\varepsilon_h}) |\w_h^{i+1}+{\G}|(\w_h^{i+1}+{\G}),\v) \big|\\
\medskip
\hspace{2cm}  \le \ds \sum_{\kappa \in \mathcal{T}_h} \| \varepsilon - \varepsilon_h \|_{L^\infty(\kappa)} \| \w_h^{i+1}+{\G}\|_{L^4(\kappa)^d} \| \nabla (\w_h^{i+1}+{\G}) \|_{L^2(\kappa)^d} \|\v \|_{L^4(\kappa)^d}\\  \medskip
\hspace{3cm} + \ds  \sum_{\kappa \in \mathcal{T}_h}  \| \beta(\varepsilon_h) - \beta_{\varepsilon_h} \|_{L^6(\kappa)}  \|\w_h^{i+1}+{\G}\|_{L^6(\kappa)^d} \|\w_h^{i+1}+{\G}\|_{L^2(\kappa)^d} \|  \v \|_{L^6(\kappa)^d}\\ \medskip
\hspace{2cm}  \le \ds \| \w_h^{i+1}+{\G}\|_{L^4(\Omega)^d} \|\v \|_{L^4(\Omega)^d} \sum_{\kappa \in \mathcal{T}_h}
\| \varepsilon - \varepsilon_h \|_{L^\infty(\kappa)} \| \nabla (\w_h^{i+1}+{\G}) \|_{L^2(\kappa)^d}\\ \medskip
\hspace{3cm} + \ds \|  \v \|_{L^6(\Omega)^d}  \|\w_h^{i+1}+{\G}\|_{L^6(\Omega)^d} \sum_{\kappa \in \mathcal{T}_h}  \| \beta(\varepsilon_h) - \beta_{\varepsilon_h} \|_{L^6(\kappa)}  \|\w_h^{i+1}+{\G}\|_{L^2(\kappa)^d} \\ \medskip
\hspace{2cm}  \le \ds c_4   \Big(  \sum_{\kappa\in \mathcal{T}_h} \| \varepsilon - \varepsilon_h \|^2_{L^\infty(\kappa)}  +   \| \beta(\varepsilon_h) - \beta_{\varepsilon_h} \|^2_{L^6(\kappa)}   \Big)^{1/2}.
 \end{array}
\]
%
%The fourth term can be treated as the third one and the fifth, sixth and seventh terms can be treated as the second one as $\w_h^{i}$ and $\w_h^{i+1}$ are bounded in $X_0$.
As $\p_h^{i+1}$ is bounded in $L^2(\Omega)$, the last term satisfies the following inequality:
\[
\begin{array}{rcl}
\medskip
\medskip
\big|(b_h - b)(\v, p_h^{i+1}) \big| &=& \ds \big| \int_\Omega  \div((\varepsilon_h - \varepsilon)\v) p_h^{i+1} d\x \big|\\
\medskip
&=& \ds \big| \int_\Omega
\nabla(\varepsilon_h - \varepsilon) \v \,  p_h^{i+1} d\x + \int_\Omega  (\varepsilon_h - \varepsilon) \div(\v) \,  p_h^{i+1} d\x\big|\\
\medskip
&\le& \ds c_4  \sum_{\kappa \in \mathcal{T}_h}  \big(  \|\nabla (\varepsilon_h - \varepsilon) \|_{L^3(\kappa)^d} |\v |_{L^6(\kappa)^d}+ \|\varepsilon_h - \varepsilon\|_{L^\infty(\kappa)} \|\div(\v) \|_{L^2(\kappa)}  \big) \| p_h^{i+1} \|_{L^2(\kappa)} \\
\medskip
&\le& \ds c_5 \Big[ \Big( \sum_{\kappa \in \mathcal{T}_h}   \|\nabla (\varepsilon_h -  \varepsilon) \|^2_{L^3(\kappa)^d} \Big)^{1/2} +  \Big( \sum_{\kappa \in \mathcal{T}_h}  \|\varepsilon_h - \varepsilon\|^2_{L^\infty(\kappa)}   \Big)^{1/2}   \Big]  |\v |_X.
\end{array}
 \]
 All the above inequalities leads to the following one:
\begin{equation}\label{upp1}
\begin{array}{rcl}
\ds \underset { \v \ne 0 } {\underset { \v \in X_0} {\sup} } \frac{{\mathcal F}_\varepsilon (\v)}{|\v |_{X_0}} &\le& \ds c_6 \Big(
\sum_{\kappa \in \mathcal{T}_h} \big( \| \varepsilon - \varepsilon_h \|^2_{L^\infty(\kappa)} +  \|\nabla \varepsilon_h - \nabla \varepsilon \|^2_{L^3(\kappa)^d} \\
\medskip
&& \qquad \qquad + \| \alpha (\varepsilon) - \alpha_{\varepsilon_h}\|^2_{L^2(\kappa)} + \| \beta (\varepsilon) - \beta_{\varepsilon_h}\|^2_{L^6(\kappa)}  \big)
 \Big)^{1/2}.
\end{array}
\end{equation}
We treat now the second term of the right hand side of Equation \eqref{majorationerreur}. We have by using \eqref{equatt2},
\[
\begin{array}{rcl}
\medskip
\mathcal{F}_{it} (\v) &=& -c_h(\w^{i+1}_h+{\G},\w_h^{i+1}+{\G},\v) + c_h(\w^i_h+{\G},\w_h^{i+1}+{\G},\v)\\
\medskip
&=& d_h(\w_h^i - \w^{i+1}_h,\w_h^{i+1}+{\G},\v) + (\beta_h (|\w^i_h+{\G}| - |\w^{i+1}_h+{\G}|) (\w^i_h+{\G}), \v)\\
\medskip
&\le& c_7 \|\w_h^i - \w^{i+1}_h \|_X \|\w_h^{i+1}+{\G}|_X \|\v |_X\\
&\le& c_7 \ds  \Big(  \sum_{\kappa \in \mathcal{T}_h} (\eta_{i,\kappa}^L )^2 \Big)^{1/2} |\v |_X.
\end{array}
\]
We deduce then the above inequalities
\begin{equation}\label{upp11}
\ds \underset { \v \ne 0 } {\underset { \v \in X_0} {\sup} } \frac{{\mathcal F}_{it} (\v)}{|\v |_{X_0}} \le c_7 \ds  \Big(  \sum_{\kappa \in \mathcal{T}_h} (\eta_{i,\kappa}^L )^2 \Big)^{1/2}.
\end{equation}
Let us now bound the third term of the right hand side of \eqref{majorationerreur}. We take $\v_h = C_h \v $ and We use the fact that $\varepsilon_h$ is bounded by $1$ to bound each term of the  right hand side of $\mathcal{F}$ given by \eqref{FFepsivvh}.  The first term can be treated as follow:
\[
\begin{array}{rcl}
\medskip
\ds \big| \int_\Omega \varepsilon_h (\f - \f_h) (\v - \v_h) d\x  \big|  &\le&  c_8   \sum_{\kappa \in \mathcal{T}_h}   h_\kappa \| \f - \f_h |_{H^1(\kappa)^d} |\v |_{L^2(\kappa)^d} \\
&\le& c_8 \ds \Big(   \sum_{\kappa \in \mathcal{T}_h} h_\kappa^2 \| \f - \f_h \|^2_{L^2(\kappa)^d} \Big)^{1/2}  |\v |_X.
\end{array}
\]
By using \eqref{Chcont} and \eqref{equatt2}, the second, third and fourth terms satisfy
\[
\begin{array}{ll}
\medskip
\big| c_h(\w_h^i+{\G},{\G}-{\G}_h,\v_h-\v) + d_h({\G} - {\G}_h, \w_h^{i+1} + {\G}_h, \v_h - \v) \\
\hspace{2cm} +(\beta_h (|\w_h^i + {\G}| - |\w_h^i + {\G}_h|) (\w_h^{i+1} + {\G}_h), \v_h- \v) \big| \ds  \le c_9 \Big(  \sum_{\kappa \in \mathcal{T}_h} | {\G} - {\G}_h |^2_{H^1(\kappa)^d} \Big)^{1/2} |\v |_X.
\end{array}
\]
It is easy to check that the fifth term can be bounded as
\[
\big| a_h({\G}-{\G}_h,\v-\v_h)  \big| \le \ds c_{10} \Big(  \sum_{\kappa \in \mathcal{T}_h} | {\G} - {\G}_h |^2_{H^1(\kappa)^d} \Big)^{1/2} |\v |_X.
\]
We use Relations \eqref{clement1} to bound the last term of Equation \eqref{FFepsivvh} denoted by $T_e$ as following:
\[
\begin{array}{rcl}
| T_e | &\le& \ds c_{11} \Big(  \sum_{\kappa \in \mathcal{T}_h}   (\eta_{i,\kappa}^{D})^2 \Big)^{1/2}.
\end{array}
\]
All the above bounds corresponding to the terms of \eqref{FFepsivvh} allow us to get
\begin{equation}\label{upp2}
\ds \underset { \v \ne \0 } {\underset { \v \in X_0} {\sup} }  \frac{{\mathcal F}(\v_h - \v)}{ |(\v,q) |_{X_0 \times M}}  \le c_{11} \Big(  \sum_{\kappa \in \mathcal{T}_h}   (\eta_{i,\kappa}^{D})^2 + \sum_{\kappa \in \mathcal{T}_h} h_\kappa^2 \| \f - \f_h \|^2_{L^2(\kappa)^d} +  \sum_{\kappa \in \mathcal{T}_h} | {\G} - {\G}_h |^2_{H^1(\kappa)^d}  \Big)^{1/2}.
\end{equation}
To end the prove of the theorem, we need to bound the last term of Equation \eqref{majorationerreur}. The term $(b-b_h)(\w_h^{i+1},q)$ can be bounded as the last term of Equation \eqref{FFepsi} and we get
\[
| (b-b_h)(\w_h^{i+1},q) | \le c_{12} \Big[ \Big( \sum_{\kappa \in \mathcal{T}_h}   \|\nabla (\varepsilon_h -  \varepsilon) \|^2_{L^3(\kappa)^d} \Big)^{1/2} +  \Big( \sum_{\kappa \in \mathcal{T}_h}  \|\varepsilon_h - \varepsilon\|^2_{L^\infty(\kappa)}   \Big)^{1/2}   \Big]  \| q \|_{L^2(\Omega)}.
\]
Furthermore, the term $b_h(\w_h^{i+1},q)$ satisfies by using the properties of $\varepsilon_h$,
\[
\begin{array}{rcl}
|b_h(\w_h^{i+1},q) | &=& \ds  \big|  \sum_{\kappa \in \mathcal{T}_h} \int_\kappa \div(\varepsilon_h \w_h^{i+1}) \, q d\x\big| \\
\medskip
&\le& c_{13} \ds \Big( \sum_{\kappa \in \mathcal{T}_h}  \| \div(\varepsilon_h \w_h^{i+1}) \| _{L^2(\kappa)}  \Big)^{1/2} \| q \|_{L^2(\Omega)}.
\end{array}
\]
Then we obtain
\begin{equation}\label{upp3}
\begin{array}{ll}
\underset { q \ne 0 } {\underset { q \in M} {\sup} } \ds \frac{|(b-b_h)(\w_h^{i+1},q) + b_h(\w_h^{i+1},q)|}   {||q||_{M}}  \\
\medskip
\hspace{1cm} \le \ds \Big[  \sum_{\kappa \in \mathcal{T}_h}   \|\nabla (\varepsilon_h -  \varepsilon) \|^2_{L^3(\kappa)^d}  +   \sum_{\kappa \in \mathcal{T}_h}  \|\varepsilon_h - \varepsilon\|^2_{L^\infty(\kappa)} +  \sum_{\kappa \in \mathcal{T}_h}  \| \div(\varepsilon_h \w_h^{i+1}) \|^2 _{L^2(\kappa)}     \Big]^{1/2}.
\end{array}
\end{equation}
Thus Equations \eqref{upp1}, \eqref{upp11}, \eqref{upp2} and \eqref{upp3} end the prove of the theorem.\hfill$\Box$\\

Now, we address the efficiency of the indicators.
{\thm \label{theo2}
%We assume that the solution $\u$ of Problem $(FV)$ is such that $ \nabla \u  \in L^\infty(\Omega)^{d\times d}.$
For each $\kappa \in \mathcal T_h$,  we have the following estimations:
\begin{equation} \label{majoindicaline}
\begin{array}{ll}
(\eta^L_{i,\kappa})^2 \le 2  \| \w- \w^{i+1}_h \|^2_{H^1(\kappa)^d} + 2 \| \w- \w^i_h \|^2_{H^1(\kappa)^d},
\end{array}
\end{equation}
%----------------------------
%-----------------------------
\begin{equation} \label{theo1}
\mbox{ } \hspace{-4.5cm} ( \eta^{D}_{i,\kappa})^2  \le  C L(w_\kappa),
\end{equation}
where $w_\kappa$ denotes the set of elements $K \in \mathcal{T}_h$ such that $\kappa \cup K \ne \phi$ and where, for any set $W$, we define
\begin{equation} \label{LW}
\begin{array}{rcl}
\medskip
L(W) & =  &  \ds \sum_{\kappa_1 \in W}h^2_{\kappa_1} \| \f -\f_h \|^2_{L^2(\kappa_1)^d} + \| {\G}-{\G}_h\|^2_{H^1(W)^d} + \|{\G}-{\G}_h\|^2_{L^3(W)^d} \\
\medskip
&& + \| \varepsilon - \varepsilon_h \|^2_{L^\infty(W)} +  \|\nabla \varepsilon_h - \nabla \varepsilon \|^2_{L^3(W)^d} + \| \alpha (\varepsilon) - \alpha_{\varepsilon_h}\|^2_{L^2(W)} + \| \beta (\varepsilon) - \beta_{\varepsilon_h}\|^2_{L^6(W)}\\
&& +  \| \w-\w_h^{i+1}\|^2_{H^1(W)^d} +  \| \w-\w_h^i\|^2_{H^1(W)^d} +  \| \w-\w_h^{i+1}\|^2_{L^3(W)^d} + \|p - p_h^{i+1}\|^2_{L^2(W)}.
\end{array}
\end{equation}
}
\noindent \textbf {Proof.}
The estimation of the linearization indicator follows easily from the triangle inequality by introducing  $\u$ in  $(\eta^L_{i,\kappa})^2 $.\\
\noindent Now, we aim to bound $ ( \eta^{D}_{i,\kappa})^2 $. We proceed in three steps.\\
{\bf \underline{First step :}} We bound the first part of the indicator $\eta^D_{i,\kappa}$\\
\noindent  We take $\v_h=\0$ in \eqref{equation1}, we obtain:
\begin{equation}\label{equatindicat}
\begin{array}{ll}
\medskip
\ds \sum_{\kappa \in \mathcal{T}_h} \Big\{ \int_\kappa \Big(\varepsilon_h \f_h + \frac{1}{Re}  \div (\varepsilon_h  \nabla (\w_h^{i+1} + {\G}_h)  ) - \alpha_{\varepsilon_h} (\w_h^{i+1}+{\G}_h) - \beta_{\varepsilon_h} |{\G}_h + \w_h^i|  ({\G}_h+ \w_h^{i+1})\\
\medskip
 \qquad \qquad \ds - \varepsilon_h (({\G}_h + \w_h^i)\cdot \nabla ) ({\G}_h+\w_h^{i+1}) -
\frac{1}{2} \div(\varepsilon_h ({\G}_h + \w_h^i)) ({\G}_h+ \w_h^{i+1} ) - \varepsilon_h \nabla p_h^{i+1}
%\\
%\medskip
% \qquad \qquad
 \Big) \v  d\x\\
\medskip
 \qquad \quad - \ds \frac{1}{2} \sum_{e \in \varepsilon_{\kappa}\cap \mathcal{E}_h} \int_{e}  \ds [(\frac{1}{Re} \varepsilon_h \nabla (\w_h^{i+1} + {\G}_h) - p_h^{i+1} \mathbb I)(\sigma)\cdot {\n}] \v d \sigma
\Big\}\\
\medskip
= -( \varepsilon_h (\f - \f_h) , \v)  + c_h(\w_h^i+{\G},{\G}-{\G}_h,\v) + d_h({\G} - {\G}_h, \w_h^{i+1}+{\G}_h, \v) \\
\medskip
\; \; \; \; + (\beta_h ( |\w_h^i + {\G} | - |\w_h^{i+1} + {\G}_h |) (\w_h^{i+1}+{\G}_h),\v) +a_h({\G}-{\G}_h,\v) \\
\medskip
\; \; \; \; + a(\w -\w_h^{i+1},\v) +  c(\w+{\G},\w - \w_h^{i+1},\v)
+ d(\w-\w_h^{i+1},\w_h^{i+1}+{\G},\v) \\
\medskip
\; \; \; \; +(\beta(|\w+{\G}|-|\w_h^{i+1} + {\G}|) (\w_h^{i+1} + {\G}),\v) - b(\v,p-p_h^{i+1}) - \mathcal{F}_\varepsilon (\v) - \mathcal{F}_{it} (\v)
\end{array}
\end{equation}
\vskip.2cm
\noindent We denote by $T_i, i=1, 12$ the terms of the right hand side of the last equation.\\
\noindent For a given $\kappa \in \mathcal{T}_h$, we choose $\v=\v_\kappa$ as follows:
\begin{equation*}
\v_{\kappa}=
\left \{
\begin{array}{lcl}
\Big(\varepsilon_h \f_h + \frac{1}{Re}  \div (\varepsilon_h  \nabla (\w_h^{i+1} + {\G}_h)  ) - \alpha_{\varepsilon_h} (\w_h^{i+1}+{\G}_h) - \beta_{\varepsilon_h} |{\G}_h + \w_h^i|  ({\G}_h+ \w_h^{i+1}) &\\
\medskip
 \qquad  \ds - \varepsilon_h (({\G}_h + \w_h^i)\cdot \nabla ) ({\G}_h+\w_h^{i+1}) -
\frac{1}{2} \div(\varepsilon_h ({\G}_h + \w_h^i)) ({\G}_h+ \w_h^{i+1}) - \varepsilon_h \nabla p_h^{i+1}
%\\
%\medskip
% \qquad \qquad
 \Big) \psi_{\kappa} & \hspace{-0.5cm} \mbox{on } \kappa, \\
0 & \hspace{-0.3cm}\mbox{ on } \Omega \backslash \kappa, \\
\end{array}
 \right.
\end{equation*}
where  $\psi_{\kappa}$ is the bubble function of the element  $\kappa$. Since this function vanishes outside $\kappa$ as well as on its boundaries, all volume integrals in \eqref{equatindicat} reduce to integrals on $\kappa$, and the edge (or face) integrals vanish.\\
\noindent We obtain then
\begin{equation}\label{low1}
\begin{array}{ll}
\medskip
\|\Big(\varepsilon_h \f_h + \frac{1}{Re}  \div (\varepsilon_h  \nabla (\w_h^{i+1} + {\G}_h)  ) - \alpha_{\varepsilon_h} (\w_h^{i+1}+{\G}_h) - \beta_{\varepsilon_h} |{\G}_h + \w_h^i|  ({\G}_h+ \w_h^{i+1}) &\\
\medskip
\ds - \varepsilon_h (({\G}_h + \w_h^i)\cdot \nabla ) ({\G}_h+\w_h^{i+1}) -
\frac{1}{2} \div(\varepsilon_h ({\G}_h + \w_h^i)) ({\G}_h+ \w_h^{i+1})  - \varepsilon_h \nabla p_h^{i+1}
 \Big)  \psi^{1/2} \|_{L^2(\kappa)^d}
 \le \ds \sum_{i=1}^{12} | T_i |.
\end{array}
\end{equation}
\noindent We will bound each term of the right hand side of the last equation. For this purpose we denote by $(a,b)_\kappa= \ds \int_\kappa a \, b \, d\x$, use the Cauchy-Schawrtz inequality, the fact that ${\G} \in H^1(\Omega)^d$ and that $\w_h^{i}$ is bounded in $H^1(\Omega)^d$, and the fact that $\v_\kappa$ is in $H^1_0 (\Omega)$ and that its support is $\kappa$ to show, for $1\le p \le 6$,  that
$$|| \v_{\kappa}||_{L^p(\kappa)^d} = || \v_{\kappa}||_{L^p(\Omega)^d} \leq S_p | \v_{\kappa}|_{H^1(\Omega)^d} = S_p | \v_{\kappa}|_{H^1(\kappa)^d}.$$
We also use some times the properties of the Cl\'ement operator \eqref{clement1} to bound the term $\G_h$ in $H^1(\Omega)^d$ and the inverse inequality \eqref{dk} with this $\| \v \|_{L^6(\kappa)^d} \le C_I^0(6) |\v_h |_{H^1(\kappa)^d}$.\\

\noindent The first term satisfies the following inequality:
\[
|T_1| = \big| (\varepsilon_h (\f-\f_h),\v)_\kappa \big| \le \ds
 \| \f -\f_h \|_{L^2(\kappa)^d} \|\v \|_{L^2(\kappa)^d}.
\]
The second term can be treated as following:
\[
\begin{array}{rcl}
\medskip
T_2 &=& \ds ( \varepsilon_h ((\w_h^i+{\G})\cdot \nabla ) ({\G}-{\G}_h) , \v_\kappa)_\kappa + \ds \frac{1}{2} (\div(\w_h^i+{\G})({\G}-{\G}_h) , \v_\kappa)_\kappa \\
\medskip
&& + (\beta_h |\w_h^i+{\G}| ({\G}-{\G}_h) , \v_\kappa)_\kappa\\
\medskip
&\le& \|\w_h^i+{\G} \|_{L^3(\kappa)^d} | {\G}-{\G}_h|_{H^1(\kappa)^d} \|\v_\kappa \|_{L^6(\kappa)^d} + \ds \frac{\sqrt{d}}{2} |\w_h^i+{\G} |_{H^1(\kappa)^d} \|{\G}-{\G}_h\|_{L^3(\kappa)^d} \|\v_\kappa \|_{L^6(\kappa)^d}\\
\medskip
&& + \beta_M \|\w_h^i+{\G} \|_{L^2(\kappa)^d} \| {\G}-{\G}_h\|_{L^3(\kappa)^d} \|\v_\kappa \|_{L^6(\kappa)^d}\\
&\le& c_1 \big(| {\G}-{\G}_h|_{H^1(\kappa)^d} + \|{\G}-{\G}_h\|_{L^3(\kappa)^d} \big) h_\kappa^{-1} \|\v \|_{L^2(\kappa)^d}.
\end{array}
\]
The third and fourth terms can treated in similar way to $T_2$ and we get
\[
|T_3 + T_4 | \le c_2 \big(| {\G}-{\G}_h|_{H^1(\kappa)^d} +\| {\G}-{\G}_h\
|_{L^3(\kappa)^d} \big)   h_\kappa^{-1} \|\v \|_{L^2(\kappa)^d}.
\]
By using \eqref{eq:ineg_verfurth}, the fifth and sixth terms can be easily bounded as
\[
\begin{array}{rcl}
\medskip
|T_5 + T_6 | &\le& c_3 (\| {\G}-{\G}_h\|_{H^1(\kappa)^d} + \| \w-\w_h^{i+1}\|_{H^1(\kappa)^d}) |\v |_{H^1(\kappa)^d}\\
&\le& c_4 (\| {\G}-{\G}_h\|_{H^1(\kappa)^d} +  \| \w-\w_h^{i+1}\|_{H^1(\kappa)^d}) h_\kappa^{-1} \|\v \|_{L^2(\kappa)^d}
\end{array}
\]
Let us now bound the terms $T_7,T_8$ and $T_9$. By using the Green formula We have:
\[
\begin{array}{ll}
\medskip
T_7 + T_8 + T_9 = \ds \frac{1}{2} \big( \varepsilon ((\w +{\G})\cdot \nabla (\w - \w_h^{i+1}),\v_\kappa)_\kappa -  \frac{1}{2} (\varepsilon (\w +{\G})\cdot \nabla \v_\kappa, \w - \w_h^{i+1})_\kappa \big) \\
\medskip
\hspace{3cm}+ (\beta_h |\w +{\G}| (\w - \w_h^{i+1}),\v_\kappa)_\kappa + \ds \frac{1}{2} (\varepsilon(\w-\w_h^{i+1})\cdot  (\w_h^{i+1}+{\G}), \v_\kappa )_\kappa \\
\hspace{3cm} \ds -  \frac{1}{2}  (\varepsilon (\w-\w_h^{i+1}) \cdot \nabla \v_\kappa, \w_h^{i+1}+{\G})_\kappa + \ds (\beta (|\w+{\G}|-|\w_h^{i+1} + {\G}|) (\w_h^{i+1} + {\G}),\v_\kappa)_\kappa
\end{array}
\]
Then we get
\[
\begin{array}{ll}
\medskip
T_7 + T_8 + T_9 \le \\
\medskip
\hspace{.2cm} \ds \frac{1}{2} \|{\G} + \w \|_{L^3(\kappa)^d} |\w-\w_h^{i+1} |_{H^1(\kappa)^d} \|\v_\kappa \|_{L^6(\kappa)^d} + \frac{1}{2} \|{\G} + \w \|_{L^6(\kappa)^d} |\w-\w_h^{i+1} |_{H^1(\kappa)^d} \| \v_\kappa\|_{L^3(\kappa)^d}\\
\medskip
\hspace{.2cm} \ds +\beta_M \|{\G} + \w \|_{L^3(\kappa)^d} \|\w-\w_h^{i+1} \|_{L^2(\kappa)^d} \|\v_\kappa \|_{L^6(\kappa)^d} + \frac{1}{2} \|{\G} + \w_h^{i+1} \|_{L^3(\kappa)^d} \|\w-\w_h^{i+1} \|_{L^2(\kappa)^d} \| \v_\kappa\|_{L^6(\kappa)^d}\\
\medskip
\hspace{.2cm} \ds  +\frac{1}{2} \|\w-\w_h^{i+1} \|_{L^3(\kappa)^d} |\v_\kappa |_{H^1(\kappa)^d}  \|{\G} + \w_h^{i+1} \|_{L^6(\kappa)^d} +\beta_M |\w-\w_h^{i+1} |_{L^2(\kappa)^d}  \|{\G} + \w_h^{i+1} \|_{L^6(\kappa)^d} \|\v_\kappa \|_{L^3(\kappa)^d}
\end{array}
\]
and then
\[
T_7+T_8+T_9 \le c_5 \big(  \|\w-\w_h^{i+1} \|_{H^1(\kappa)^d} + \|\w-\w_h^{i+1} \|_{L^3(\kappa)^d}  \big)  h_\kappa^{-1} |\v_\kappa |_{L^2(\kappa)^d}
\]
The tenth term $T_{10}$ can be bounded as follows:
\[
\begin{array}{rcl}
\medskip
T_{10} &=& (\div(\varepsilon \v_\kappa)_\kappa (p-p_h^{i+1}))_\kappa\\
\medskip
&\le& (\varepsilon \div(\v_\kappa) (p-p_h^{i+1}))_\kappa + (\nabla \varepsilon \cdot \v_\kappa, (p-p_h^{i+1}))_\kappa\\
\medskip
&\le& c_6 (\|\varepsilon \|_{L^\infty(\kappa)} |\v_\kappa |_{H^1(\kappa)^d} + \|\nabla \varepsilon \|_{L^3(\kappa)} |\v_\kappa |_{L^6(\kappa)^d} ) \|p - p_h^{i+1}\|_{L^2(\kappa)}\\
&\le& c_7 h_\kappa^{-1} |\v_\kappa |_{L^2(\kappa)^d} ) \|p - p_h^{i+1}\|_{L^2(\kappa)}.
\end{array}
\]
It remains to bound the last two term of \eqref{equatindicat}. \\
\noindent We begin with $\mathcal{F}_{it}(\v)$:
\[
\begin{array}{rcl}
\medskip
\mathcal{F}_{it}(\v) &=& -c_h(\w^{i+1}_h+{\G},\w_h^{i+1}+{\G},\v_\kappa) + c_h(\w^i_h+{\G},\w_h^{i+1}+{\G},\v_\kappa)\\
\medskip
&=& d_h(\w_h^i - \w^{i+1}_h,\w_h^{i+1}+{\G},\v_\kappa) + (\beta_h (|\w^i_h+{\G}| - |\w^{i+1}_h+{\G}|) (\w^i_h+{\G}), \v_\kappa)\\
\medskip
&\le& c_7 \|\w_h^i - \w^{i+1}_h \|_{H^1(\kappa)} \|\w_h^{i+1}+{\G}\|_{H^1(\kappa)} |\v_\kappa |_{H^1(\kappa)}\\
&\le& c_8 \ds   \eta_{i,\kappa}^L   h_\kappa^{-1} |\v_\kappa |_{L^2(\kappa)^d}.
\end{array}
\]
The term $\mathcal{F}_\varepsilon$ can be bounded exactly as in the proof of Theorem \ref{upperb} (see \eqref{upp1}) and we get
\[
\mathcal{F}_\varepsilon \le \ds c_9
 \big( \| \varepsilon - \varepsilon_h \|^2_{L^\infty(\kappa)} +  \|\nabla \varepsilon_h - \nabla \varepsilon \|^2_{L^3(\kappa)^d} + \| \alpha (\varepsilon) - \alpha_{\varepsilon_h}\|^2_{L^2(\kappa)} + \| \beta (\varepsilon) - \beta_{\varepsilon_h}\|^2_{L^6(\kappa)}\big) h_\kappa^{-1} |\v_\kappa |_{L^2(\kappa)^d}.
\]
Finally, by regrouping all the above inequalities, Relation \eqref{low1} multiplied by $h_\kappa$ gives by using Property \ref{psi},
\begin{equation}\label{ubound1}
\begin{array}{ll}
\medskip
h_\kappa \|\varepsilon_h \f_h + \frac{1}{Re}  \div (\varepsilon_h  \nabla (\w_h^{i+1} + {\G}_h)  ) - \alpha_{\varepsilon_h} (\w_h^{i+1}+{\G}_h) - \beta_{\varepsilon_h} |{\G}_h + \w_h^i|  ({\G}_h+ \w_h^{i+1}) &\\
\medskip
 \qquad  \ds - \varepsilon_h (({\G}_h + \w_h^i)\cdot \nabla ) ({\G}_h+\w_h^{i+1}) -
\frac{1}{2} \div(\varepsilon_h ({\G}_h + \w_h^i)) ({\G}_h+ \w_h^{i+1}) - \varepsilon_h \nabla p_h^{i+1}\|_{L^2(\kappa)^d}  \\
 \medskip
 \le c \Big( h_\kappa \| \f -\f_h \|_{L^2(\kappa)^d} + |{\G}-{\G}_h|_{H^1(\kappa)^d} + \|{\G}-{\G}_h\|_{L^3(\kappa)^d} \\
 \qquad \qquad +  \| \varepsilon - \varepsilon_h \|_{L^\infty(\kappa)} +  \|\nabla \varepsilon_h - \nabla \varepsilon \|_{L^3(\kappa)^d} + \| \alpha (\varepsilon) - \alpha_{\varepsilon_h}\|_{L^2(\kappa)} + \| \beta (\varepsilon) - \beta_{\varepsilon_h}\|_{L^6(\kappa)} \\
 \qquad \qquad +  \| \w-\w_h^{i+1}\|_{H^1(\kappa)^d} +  \| \w-\w_h^{i+1}\|_{L^3(\kappa)^d} +  \|\w-\w_h^i\|_{H^1(\kappa)^d} + \|p - p_h^{i+1}\|_{L^2(\kappa)}  \Big).
\end{array}
\end{equation}
{\bf \underline{Second step :}} We bound the second part of the indicator $\eta^D_{i,\kappa}$.\\
\noindent Rewriting Equation \eqref{equatindicat} we infer
\begin{equation}\label{equatindicatt}
\begin{array}{ll}
\medskip
\ds \frac{1}{2} \sum_{\kappa \in \mathcal{T}_h} \sum_{e \in \varepsilon_{\kappa}\cap \mathcal{E}_h} \int_{e}  \ds [(\frac{1}{Re} \varepsilon_h \nabla (\w_h^{i+1} + {\G}_h) - p_h^{i+1} \mathbb I)(\sigma)\cdot {\n}] \v d \sigma \\
\medskip
= \ds \sum_{\kappa \in \mathcal{T}_h} \int_\kappa \Big(\varepsilon_h \f_h + \frac{1}{Re}  \div (\varepsilon_h  \nabla (\w_h^{i+1} + {\G}_h)  ) - \alpha_{\varepsilon_h} (\w_h^{i+1}+{\G}_h) - \beta_{\varepsilon_h} |{\G}_h + \w_h^i|  ({\G}_h+ \w_h^{i+1})\\
\medskip
 \qquad \qquad \ds - \varepsilon_h (({\G}_h + \w_h^i)\cdot \nabla ) ({\G}_h+\w_h^{i+1}) -
\frac{1}{2} \div(\varepsilon_h ({\G}_h + \w_h^i)) ({\G}_h+ \w_h^{i+1}) - \varepsilon_h \nabla p_h^{i+1}
 \Big) \v  d\x\\
\medskip
\; \; \; \; + ( \varepsilon_h (\f - \f_h) , \v)  - c_h(\w_h^i+{\G},{\G}-{\G}_h,\v) - d_h({\G} - {\G}_h, \w_h^{i+1}+{\G}_h, \v) \\
\medskip
\; \; \; \; - (\beta_h ( |\w_h^i + {\G} | - |\w_h^{i+1} - {\G}_h |) (\w_h^{i+1}+{\G}_h),\v) -a_h({\G}-{\G}_h,\v) \\
\medskip
\; \; \; \; - a(\w -\w_h^{i+1},\v) -  c(\w+{\G},\w - \w_h^{i+1},\v)
- d(\w-\w_h^{i+1},\w_h^{i+1}+{\G},\v) \\
\medskip
\; \; \; \; -(\beta(|\w+{\G}|-|\w_h^{i+1} + {\G}|) (\w_h^{i+1} + {\G}),\v) + b(\v,p-p_h^{i+1}) + \mathcal{F}_\varepsilon (\v) + \mathcal{F}_{it} (\v).
\end{array}
\end{equation}
\noindent For a given $e \not \subset \Gamma$, we denote by $(\kappa,\kappa') \in (\mathcal{T}_h)^2$ the two elements that share $e$ and set
\begin{equation*}
\v=\v_{e}=
\left \{
\begin{array}{lcl}
\mathcal{L}_{e} \Big(
 [\frac{1}{Re} \varepsilon_h \nabla (\w_h^{i+1} + {\G}_h) - \bar p_h^{i} \mathbb I) {\n}] \psi_{e}  \Big)& \hspace{-0.5cm} \mbox{on } \{ \kappa, \kappa' \}, \\
0 & \hspace{-0.3cm}\mbox{ on } \Omega \backslash \{ \kappa \bigcup \kappa' \}, \\
\end{array}
 \right.
\end{equation*}
\noindent where $\psi_e$ is the edge-bubble (or face-bubble) function and $\mathcal{L}_{e}$ the lifting operator.\\
We consider Equation \eqref{equatindicatt} and we replace $\v=\v_e$ to get
\begin{equation}\label{equatindicattt}
\begin{array}{ll}
\medskip
\ds \frac{1}{2} \ds \Big|\Big| [(\frac{1}{Re} \varepsilon_h \nabla (\w_h^{i+1} + {\G}_h) - p_h^{i+1} \mathbb I)\cdot {\n}] \psi^{1/2}_e  \Big|\Big|^2_{L^2(e)^d} \\
\medskip
= \ds \int_{\kappa  \cup \kappa'} \Big(\varepsilon_h \f_h + \frac{1}{Re}  \div (\varepsilon_h  \nabla (\w_h^{i+1} + {\G}_h)  ) - \alpha_{\varepsilon_h} (\w_h^{i+1}+{\G}_h) - \beta_{\varepsilon_h} |{\G}_h + \w_h^i|  ({\G}_h+ \w_h^{i+1})\\
\medskip
 \qquad  \ds - \varepsilon_h (({\G}_h + \w_h^i)\cdot \nabla ) ({\G}_h+\w_h^{i+1}) -
\frac{1}{2} \div(\varepsilon_h ({\G}_h + \w_h^i)) ({\G}_h+ \w_h^{i+1}) - \varepsilon_h \nabla p_h^{i+1}
 \Big) \v_e  d\x\\
\medskip
\quad + ( \varepsilon_h (\f - \f_h) , \v_e)  - c_h(\w_h^i+{\G},{\G}-{\G}_h,\v_e) - d_h({\G} - {\G}_h, \w_h^{i+1}+{\G}_h, \v_e) \\
\medskip
\; \; \; \; - (\beta_h ( |\w_h^i + {\G} | - |\w_h^{i+1} - {\G}_h |) (\w_h^{i+1}+{\G}_h),\v_e) -a_h({\G}-{\G}_h,\v_e) \\
\medskip
\; \; \; \; - a(\w -\w_h^{i+1},\v_e) -  c(\w+{\G},\w - \w_h^{i+1},\v_e)
- d(\w-\w_h^{i+1},\w_h^{i+1}+{\G},\v_e) \\
\medskip
\; \; \; \; -(\beta(|\w+{\G}|-|\w_h^{i+1} + {\G}|) (\w_h^{i+1} + {\G}),\v_e) + b(\v,p-p_h^{i+1}) + \mathcal{F}_\varepsilon (\v_e) + \mathcal{F}_{it} (\v_e).
\end{array}
\end{equation}
We use the Cauchy-Schwartz inequality, multiply by $h_e$, bound all the terms of the second member exactly as we did in the previous step,  use Property \ref{psii}, simplify by $h_e^{1/2} || \v_{e} ||_{L^2(e)^d }$, square the resulting inequality and, for a given $\kappa$, we sum over $ e \in \partial \kappa$; we obtain:
\begin{equation}\label{ubound2}
\ds \sum_{e\in \varepsilon_\kappa \cap \mathcal{E}_h} h_e \Big|\Big| [(\frac{1}{Re} \varepsilon_h \nabla (\w_h^{i+1} + {\G}_h) - p_h^{i+1} \mathbb I)(\sigma)\cdot {\n}] \psi^{1/2}_e  \Big|\Big|^2_{L^2(e)^d} \le L(\w_\kappa),
\end{equation}
where $w_\kappa$ denotes the set of elements of $\mathcal{T}_h$ that share at least one edge  (or face when $d=3$) with $\kappa$.\\

\noindent {\bf \underline{Third step :}}  We bound in this step the last part of the indicator $\eta_{i,k}^D$. \\
\noindent Equation \eqref{equation2} gives for all $q\in M$:
\begin{equation}\label{xxx1}
\begin{array}{ll}
\ds \int_\Omega \div(\varepsilon (\w- \w^{i+1}_h)) \, q \, d\x = \ds - \int_\Omega \div((\varepsilon - \varepsilon_h) \w_h^{i+1}) \; q \, d\x - \int_\Omega \div(\varepsilon_h \w_h^{i+1}) \; q \, d\x\\
\medskip
\hspace{2cm} = \ds - \int_\Omega (\varepsilon - \varepsilon_h) \div( \w_h^{i+1}) \; q \, d\x - \int_\Omega \nabla (\varepsilon - \varepsilon_h) \cdot \w_h^{i+1} \; q \, d\x  - \ds  \int_\Omega  \div(\varepsilon_h \w_h^{i+1}) \; q \, d\x.
\end{array}
\end{equation}
\noindent We choose for a given $\kappa \in \mathcal{T}_h$ $$ q =q_{\kappa}= \div (\varepsilon_h \w^{i+1}_h) \xi_{\kappa}. $$
\noindent where $\xi_{\kappa}$ denotes the characteristic function of $\kappa$. We obtain by using $\div(\varepsilon \w)=0$,

\begin{equation*}
\begin{array}{ll}
\medskip
|| \div  (\varepsilon  \w^{i+1}_h) ||^2_{L^2(\kappa)} = - \ds \int_\kappa \varepsilon \div(\w- \w^{i+1}_h) \, \div (\varepsilon_h \w^{i+1}_h) \, d\x - \ds \int_\kappa \nabla(\varepsilon)\cdot (\w- \w^{i+1}_h) \, \div (\varepsilon_h \w^{i+1}_h) \, d\x \\
\medskip
\hspace{3cm} \ds + \int_\Omega (\varepsilon - \varepsilon_h) \div( \w_h^{i+1}) \;\div (\varepsilon_h \w^{i+1}_h) \, d\x + \int_\Omega \nabla (\varepsilon - \varepsilon_h) \cdot \w_h^{i+1} \; \div (\varepsilon_h \w^{i+1}_h) \, d\x
\\
\medskip
\hspace{1cm} \le   \Big( \|\varepsilon \|_{L^\infty(\kappa)} \sqrt{d} \| \w-  \w^{i+1}_h\|_{H^1(\kappa )^d} + \|\nabla \varepsilon \|_{L^6(\kappa)^d}  \|\w - \w_h^{i+1} \|_{L^3(\kappa)^d}\\
\medskip
\hspace{2cm} \quad +  \|\varepsilon - \varepsilon_h \|_{L^\infty(\kappa)} \sqrt{d} \| \w^{i+1}_h\|_{H^1(\kappa )^d} + \| \nabla(\varepsilon - \varepsilon_h) \|_{L^3(\kappa)^d}\| \w^{i+1}_h\|_{L^6(\kappa )^d} \Big) || \div  (\varepsilon  \w^{i+1}_h) ||_{L^2(\kappa)}.
\end{array}
\end{equation*}
We deduce the following bound
\begin{equation}\label{ubound3}
|| \div  (\varepsilon  \w^{i+1}_h) ||_{L^2(\kappa)} \le   \big( \| \w-  \w^{i+1}_h\|_{H^1(\kappa )^d} +  \|\w - \w_h^{i+1} \|_{L^3(\kappa)^d}
+  \|\varepsilon - \varepsilon_h \|_{L^\infty(\kappa)} + \| \nabla(\varepsilon - \varepsilon_h) \|_{L^3(\kappa)^d}\|  \big).
\end{equation}
Finally, collecting \eqref{ubound1}, \eqref{ubound2} and \eqref{ubound3}, we get the final result.  \hfill$\Box$
{\rmq In the upper bound of the indicator $\eta^D_{i,\kappa}$ we have the term $\| \w-\w_h^{i+1}\|_{L^3(W)^d} $. In fact, this term does not bother the optimally of the {\it a posteriori} estimate since in the upper bound of the Theorem \ref{upperb} we bound the global error $||\w- \w^{i+1}_h||_{X}$ with the indicators up to data errors, and we have
\[
||\w- \w^{i+1}_h||_{L^3(\Omega)}\le C_3 ||\w- \w^{i+1}_h||_{X}.
\]
}
\section{Numerical results}
In this section we will show numerical investigations corresponding to the {\it a posteriori} error estimate. All the following numerical computations are performed with FreeFem++ \cite{Freefem}.
We will treat two test cases: an academic test where we know the exact solution $(\u,p)$ and the flow in packed bed reactors such as those studied \cite{BDFtheornum1, BDFtheor1,vaf,abed} (see also references inside). \\

The {\it a posteriori} error estimates between the exact and numerical solutions obtained in the previous sections will be used in this part to show numerical results based on mesh adaptation. For this objective, it is convenient to compute the following expressions, $\eta_{i}^D$ and $\eta_{i}^L$,  for the indicators,
\begin{equation}\nonumber
 \eta_{i}^D=\big(\sum _{\kappa \in
 \mathcal{T}_h }(\eta_{\kappa,i}^D)^2
\big)^{\frac{1}{2}}
\end{equation}
and
\begin{equation}\nonumber
 \eta_{i}^L=\big(\sum _{\kappa \in
 \mathcal{T}_h }((\eta_{K,i}^L)^2
\big)^{\frac{1}{2}}.
\end{equation}
These indicators are used for mesh adaptation by the adapted mesh algorithm  introduced in \cite{BDMS}. For a given mesh, the iterations are stopped following the criteria
\begin{equation}
\label{eq:error_tolerance}
\eta_{i}^L\leq \tilde{\gamma} \eta_{i}^D,
%\quad
%\mbox{and} \quad \eta_{i}^{(D)} \leq %\upsilon,
\end{equation}
where $\tilde{\gamma}=0.01$. For the study of the dependence of the stopping criteria $\eqref{eq:error_tolerance}$ with $\tilde{\gamma}$, we refer to \cite{LAM11} and
\cite{ERN11} where the authors introduce this new stopping criterion and choose in practice $\tilde{\gamma} = 0.1$ in their numerical experiments. In this work, we choose the same value taken in \cite{BDMS}, $\tilde{\gamma} = 0.01$ for our numerical applications. For the adaptive mesh (refinement and coarsening), we use routines in FreeFem++. \\

At each refinement level mesh, we introduce the relative total error indicator given by
\begin{equation}\label{errtot}
E_{total}=\ds \frac{\eta^D}{|\u|_X + \|p \|_M},
\end{equation}
where $\eta^D$ is the discretization indicators after convergence on the iterations i (by using the stopping criteria \eqref{eq:error_tolerance}).
\subsection{First test case}
This academic case deals with numerical tests for a given exact solution $(\u,p)$ of Problem \eqref{V1}  where $\Omega=]0,1[^2$ and $\g=0$. we consider the following porosity
\[
\varepsilon(x,y) = \ds \frac{1+ e^{x+y}}{10},
\]
and the following Darcy and Forchheimer terms
\[
\alpha(\varepsilon) = (1-\varepsilon)^2, \quad \beta(\varepsilon) = 1+\varepsilon.
\]
The exact solution is given by $(\u,p)=(\ds \frac{1}{\varepsilon} \curl \psi, p)$ where,
\begin{equation*}
\left\{
\begin{array}{lcl}
\medskip
\psi(x,y)&=&\ds  e^{-30[(x-0.5)^2+(y-0.5)^2]},   \\
 p(x,y)&=&\cos(\pi x)\cos(\pi y).
\end{array}
\right.
\end{equation*}
We can check easily that $\u|_\Gamma = \0$ and $\div(\varepsilon \u)=0$. Furthermore, as $\g=\0$, then we can consider $\G=\0$.\\

We begin by testing the convergence of the iterative scheme \eqref{Whi} for a uniform mesh with respect to the Reynolds number $Re$. in this case, we use the following classical stopping criterion:
\[
\eta_i^L \le 1e-6.
\]
\noindent  The first numerical simulations corresponding to scheme \eqref{Whi} with $N=40$  show that the algorithm converges for $Re \le 130$. These results are coherent with the literature (see for instance \cite{BDMS} for more details) which announces that the convergence depends on the Reynolds number.\\

In  \cite{BDMS}, the authors proposed a simple modification (a relaxed numerical scheme) of their numerical scheme allowing to get convergence for a larger range of Reynolds numbers; by conducting their idea to our problem, we introduce the following relaxed iterative scheme ($\bar \w_h^i=\bar \u_h^i$ as $\G=0$):
\begin{equation}\label{Whir}
\left\{
\begin{array}{ll}
\medskip
\forall \v_h \in X_{0h}, \quad \ds   a_h(\bar \u_h^{i+1},\v_h) + d_h(\tilde{\u}_h^i,\bar \u_h^{i+1},\v_h) + (\beta (\varepsilon_h) | \bar{\u}_h^i| \bar{\u}_h^{i+1} , \v_h) - b_h(\v_h,\bar p_h^{i+1})  = \ell_h(\v_h),\\
\forall q_h\in M_{0h}, \quad \ds b_h(\bar \u_h^{i+1}, q_h) = 0,
\end{array}
\right.
\end{equation}
where
$$
\tilde  \u^{i}_h = \ds \frac{ \bar \u_h^{i}+ \tilde \u^{i-1}_h}{2}.
$$
\noindent In our tests, Scheme \eqref{Whir} for $N=40$ and with the stopping criterion $\bar \eta^L_i \le 10^{-6}$, %\eqref{critere1},
converges for large values of $Re\le 2000$. Thus, we adopt this relaxed method for numerical tests. \\

The indicators corresponding to the iterative scheme \eqref{Whir} will reduce in this case to ($\bar \u_h^{i+1}=\bar \w_h^{i+1}$)
\begin{equation}
\begin{array}{rcl}
\medskip
\eta_{i,\kappa}^L &=&   \|\bar \u_h^{i+1}- \bar \u_h^{i} \|_{H^1(\kappa)^d},\\
\medskip
\eta_{i,\kappa}^{D} &=& \ds h_{\kappa } || \varepsilon_h \f_h + \frac{1}{Re}  \div (\varepsilon_h  \nabla \bar \u_h^{i+1} ) - \alpha_{\varepsilon_h} \bar \u_h^{i+1} \ds - \varepsilon_h ( \tilde \u_h^i\cdot \nabla ) \bar \u_h^{i+1} -
\frac{1}{2} \div(\varepsilon_h  \tilde \u_h^i) \bar \u_h^{i+1} \\
\medskip%
&& \qquad \qquad - \beta_{\varepsilon_h}  |\bar \u_h^i|  \bar \u_h^{i+1}  - \varepsilon_h \nabla \bar p_h^{i+1} ||_{L^2(\kappa)^d}\\
\medskip
&& + \ds\frac{1}{2} \sum_{e \in \varepsilon_{\kappa}} h^{1/2}_{e} ||[(\frac{1}{Re} \varepsilon_h \nabla \bar \u_h^{i+1} - \bar p_h^{i+1} \mathbb I)(\sigma)\cdot {\n}] ||_{L^2(e)^d}
 + ||\div (\varepsilon_h \bar \u^{i+1}_h)||_{L^2(\kappa)},
\end{array}
\end{equation}
where $\varepsilon_h$ is the $\P_1$ finite element approximation of $\varepsilon$, that is $\varepsilon_h$ defined as the finite element interpolate of the porosity.\\
We define also the following relative total error between the exact and numerical solutions
\[
err = \ds \frac{|\u - \bar \u_h|_X+ \|p-\bar \p_h \|_M}{|\u |_X + \|p \|_M}
\]
where $(\bar \u_h,\bar p_h)$ is the numerical solution on a given mesh  and after convergence on the iterations $i$ (by using the stopping criteria \eqref{eq:error_tolerance}).\\
All computations start on a uniform initial triangular mesh obtained by dividing $\Omega$ into $N^2$ triangles (each edge is divided into $N$ segments of equal length). Furthermore, the numerical results showed in the following are performed for $\u_h^0=\0$ and $Re=500$.\\
%solution of the Stokes problem corresponding to \eqref{V1} by canceling the corresponding non-linear terms. \\

\begin{figure}[htbp]
\begin{minipage}[b]{0.40\linewidth}
 \centering
\includegraphics[width=7cm]{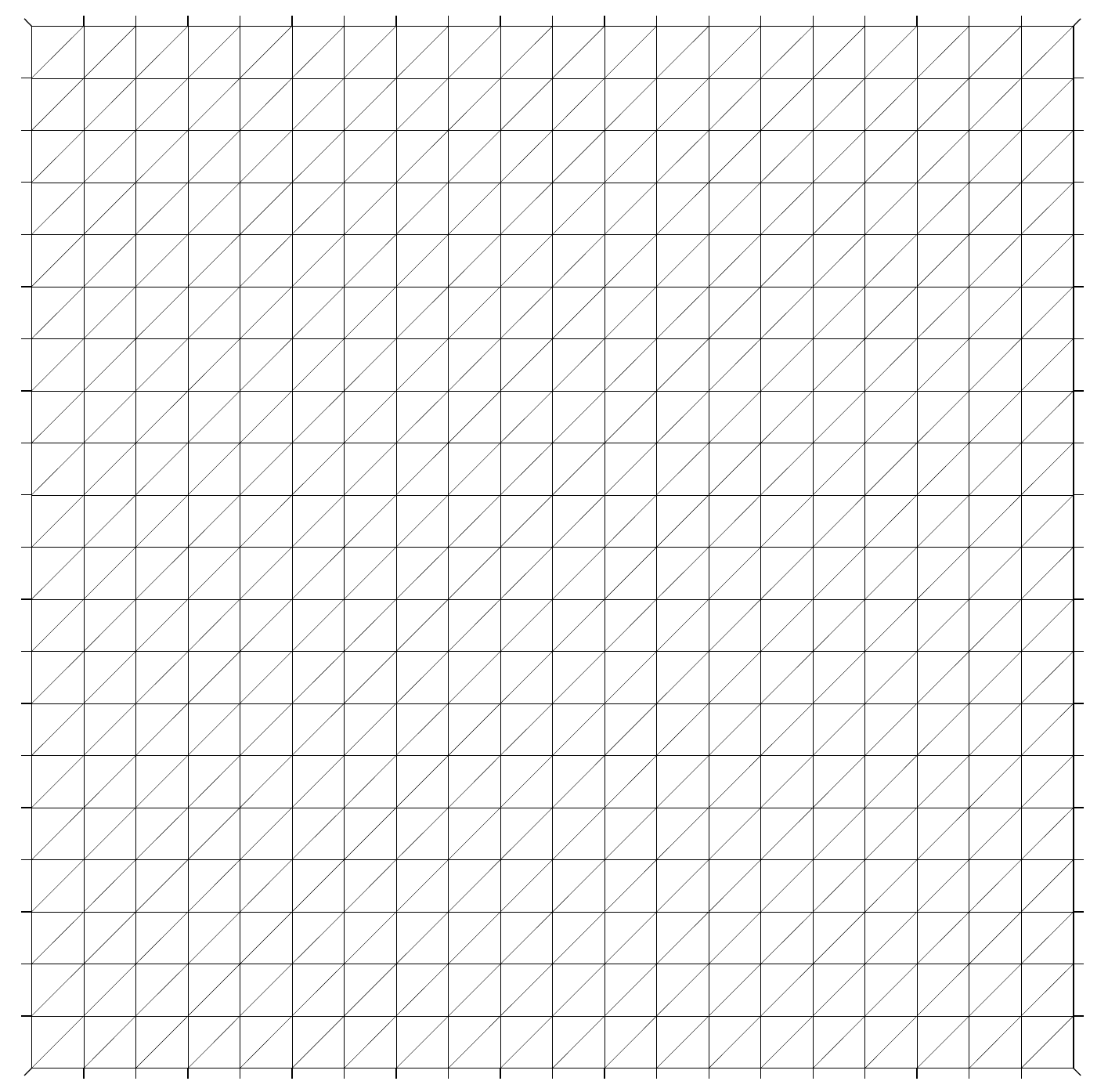}
\vspace{-.8cm}
\caption{Initial mesh (800 triangles)} \label{ figure1 }
\end{minipage}\hfill
\begin{minipage}[b]{0.40\linewidth}
\centering
\includegraphics[width=7cm]{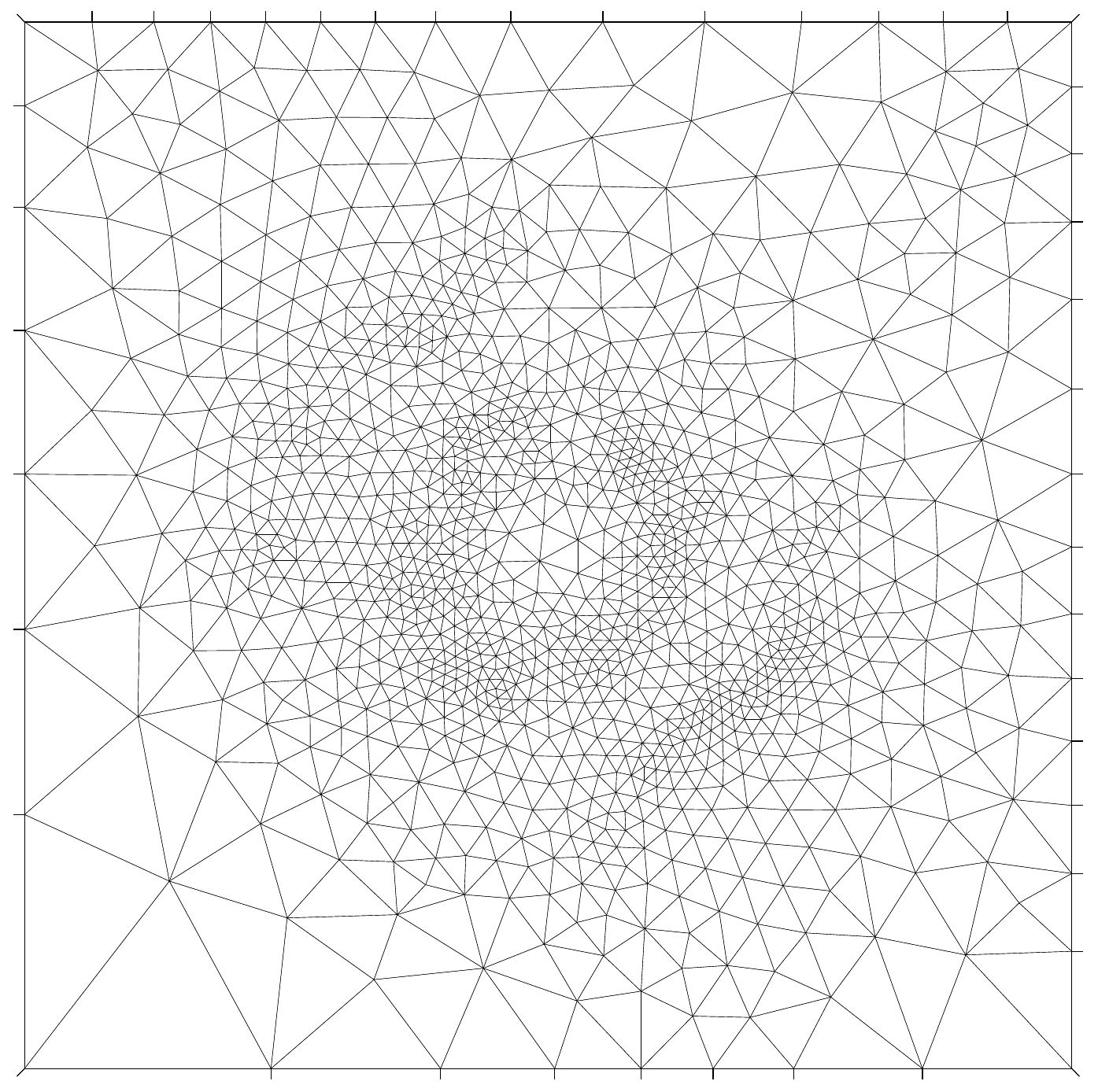}
\vspace{-.8cm}
\caption{ Second level mesh (2152 triangles)} \label{ figure2 }
\end{minipage}\hfill
\begin{minipage}[b]{0.40\linewidth}
\centering
\includegraphics[width=7cm]{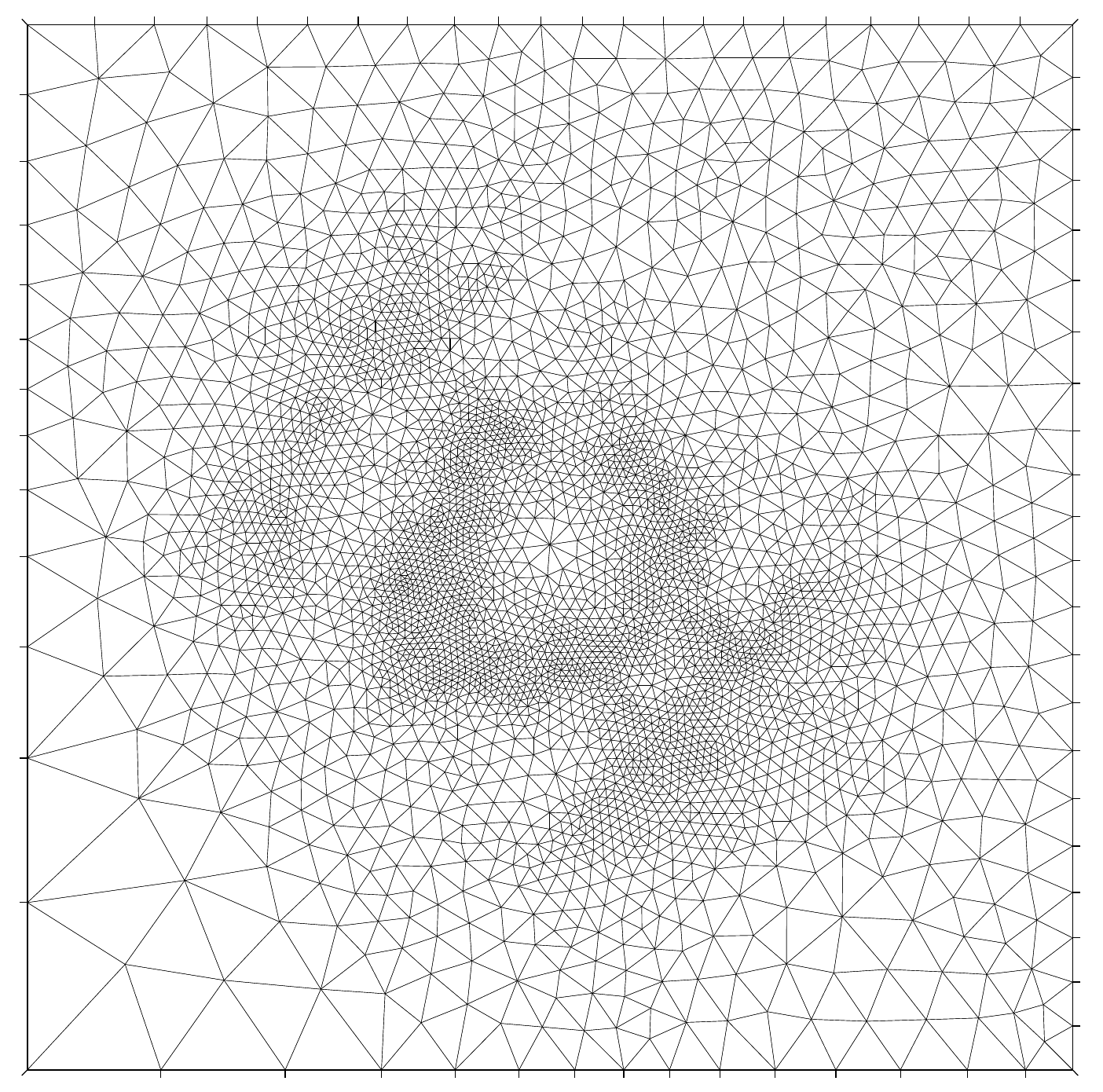}
\vspace{-.8cm}
\caption{ Fourth level mesh (6662 triangles) } \label{figure3}
\end{minipage}\hfill
\begin{minipage}[b]{0.40\linewidth}
\centering
\includegraphics[width=7cm]{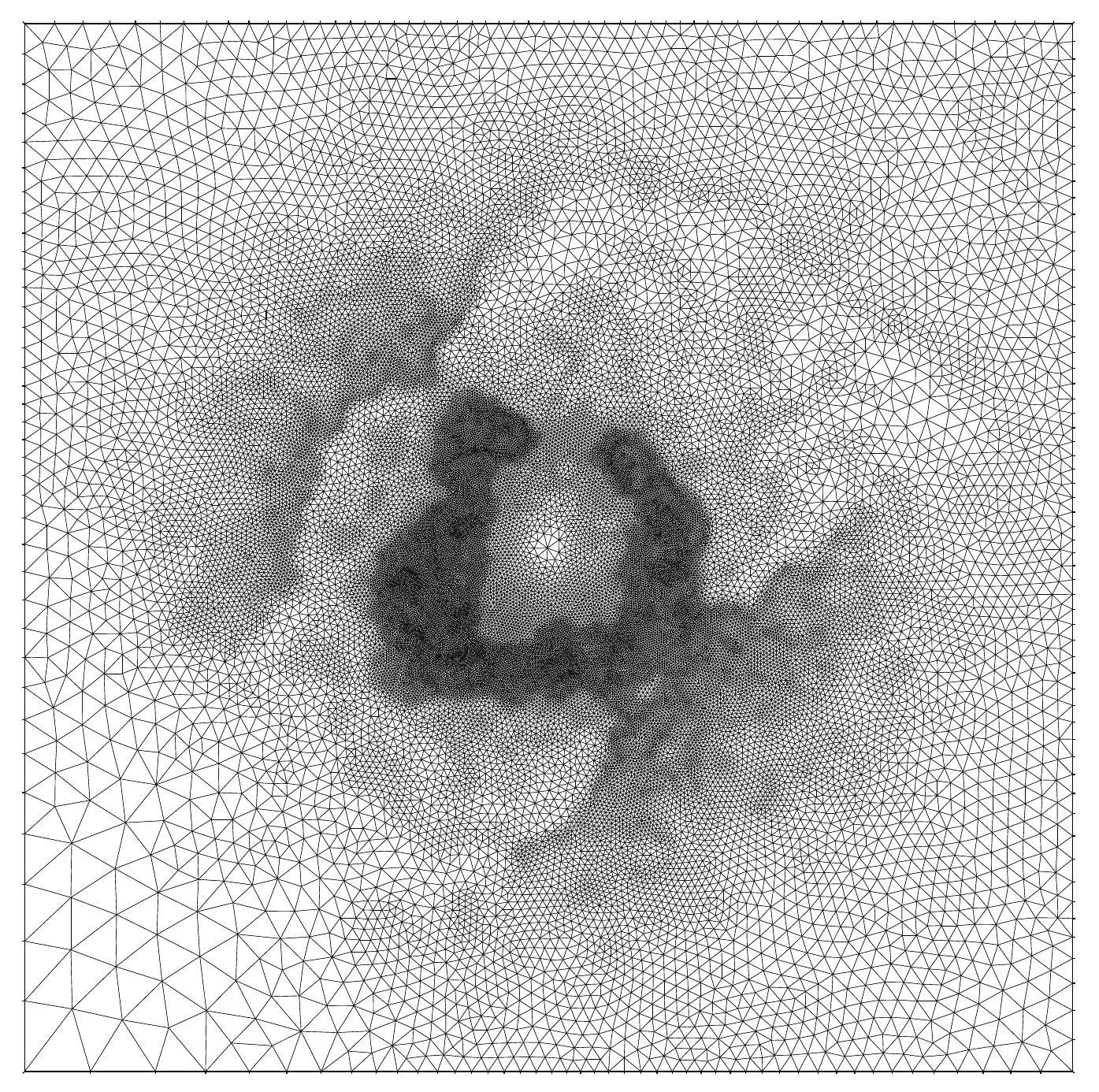}
\vspace{-.8cm}
\caption{ Eighth level mesh (61255 triangles) } \label{figure4}
\end{minipage}\hfill
\end{figure}

Figures \ref{ figure1 }-\ref{figure4} show the evolution of the mesh during the refinement levels of the algorithm when the algorithm starts with a uniform mesh produced with $N=20$.
We remark that, from an iteration to another, the concentration of the refinement is on the complex vorticity region.

%-------------------------------
\begin{figure}[!ht]
\vspace{0cm}
\begin{center}
\includegraphics[width=7.8cm]{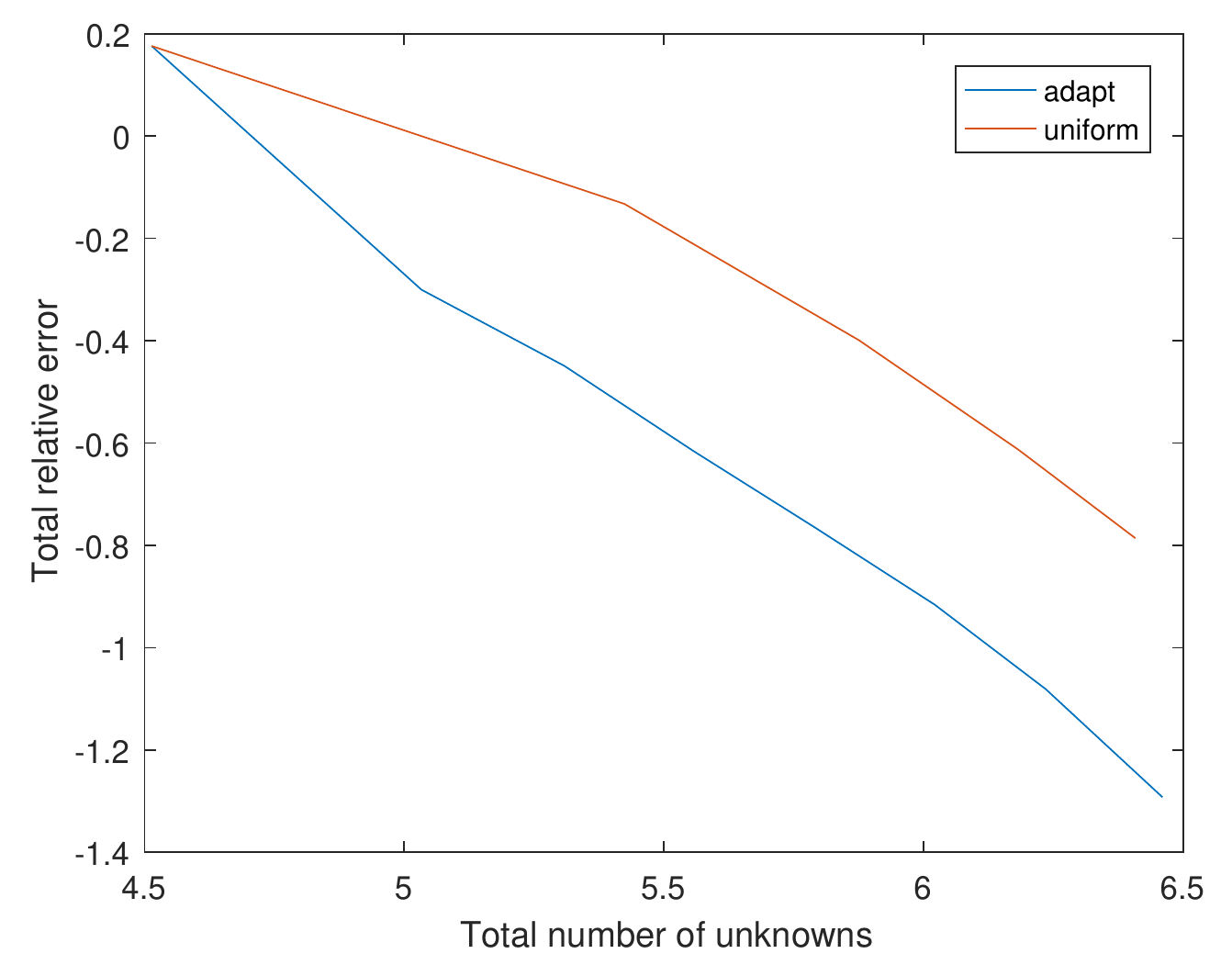}
%\hspace{0cm} \hspace{-.6cm}\caption{} %\label{globalerror}
%\end{center}
%\end{figure}
%
%\begin{figure}[!ht]
%\vspace{0cm}
%\begin{center}
\includegraphics[width=7.8cm]{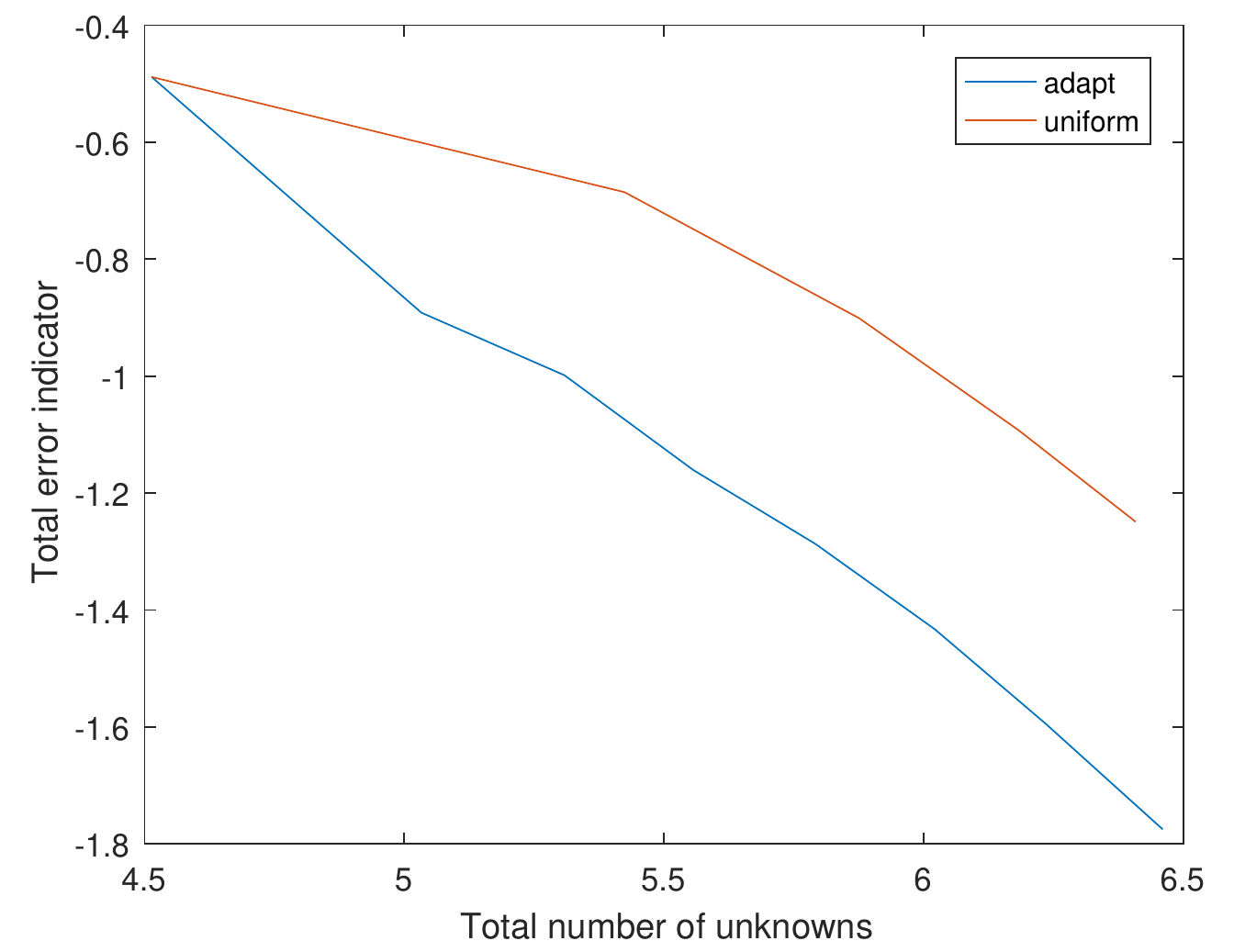}
\vspace{-0.4cm} \caption{left: Total relative errors $err$. Right: Total  error indicators $E_{total}$.} \label{Error}
\end{center}
\end{figure}

Next, in the left part of Figure \ref{Error}, we compare the relative total error $err$ with respect to the total number of unknowns in logarithmic scale, for both uniform and adaptive numerical algorithms. In the right part of Figure \ref{Error}, we show a comparison of the global error for the indicators $E_{total}$ with respect to the total number of unknowns in logarithmic scale for the adaptive and  refinement methods. Both parts of Figure \ref{Error} show clearly the advantage of the adaptive method versus the uniform one since the total errors are smaller.

We define the efficiency index by
\begin{equation}
    EI=\frac{\ds \eta^D}{\ds |\u_h-\u|_X + \|p_h - p\|_M}.
\end{equation}
In Table \ref{tab1}, we can see the value of the efficiency index for different values of unknowns STU:
\begin{center}
\begin{tabular}{|c|c|c|c|c|c|c|c|c|c|c|c|}%
\hline
STU &  32634 & 107910 & 203523 & 360143 & 619479 & 1051195 & 1720257 & 2886564 & 4853220 & 7270008  \\
\hline
EI &  0.216  & 0.256  & 0.282  & 0.284  & 0.301  & 0.304 & 0.306  & 0.328  & 0.347  & 0.371  \\
\hline
\end{tabular}
\vskip .2cm
{TABLE 1. Repartition of EI with respect to the total number of unknowns STU.}\label{tab1}
\end{center}
We can see that the coefficient EI increases from $0.216$ to $0.371$ when the STU number increases from $32634$ to $7270008$.
\subsection{Second test case}
In this section we will show numerical investigations corresponding to the flow in packed bed reactors such as those studied \cite{BDFtheornum1, BDFtheor1,vaf,abed}. In \cite{BDFtheornum1}, the authors treat a numerical iterative scheme and  study the dependency of the convergence with respect to the velocity in the inlet part of the domain for different values of the Reynolds number $Re$.

In this case, the considered smooth porosity is the following:
%(see \cite{BDFtheornum1,vaf,abed}):
\[
\varepsilon(x,y) = \ds 0.45 \big( 1 +\frac{1-0.45}{0.45} e^{-(1-y)} \big).
\]
The Darcy and Forchheimer terms are defined in \cite{abed, BDF3}:
\[
\alpha(\varepsilon) = \ds \frac{150}{Re} \big( \frac{1-\varepsilon}{\varepsilon} \big)^2, \qquad \beta(\varepsilon) = 1.75 \ds \big( \frac{1-\varepsilon}{\varepsilon} \big).
\]
It is easy to see that they satisfy all the assumptions \ref{assumpdata}.\\

\noindent We set $\Omega=]0,2[\times ]0,1[$. We distinguish between the inlet, outlet, and membrane parts of the boundary $\Gamma$ and denote them by
$\Gamma_{in}$, $\Gamma_{out}$, and $\Gamma_w$, respectively. Let
\[
\Gamma_{in} = \{0\}\times [0, 1],\quad  \Gamma_{out} = \{2\}\times [0, 1] \quad \mbox{and} \quad  \Gamma_w = [0, 2] \times (\{0\} \cup \{1\}).
\]
At the inlet $\Gamma_{in}$ and at the membrane wall $\Gamma_w$ we prescribe
Dirichlet boundary conditions, namely, the flow conditions
\[
\u|_{\Gamma_{in}} = \u_{in} =(C_{in} (1-y)y,0) \qquad \mbox{and} \qquad \u|_{\Gamma_{w}} = \0,
\]
where $C_{in} >0$. At the outlet $\Gamma_{out}$ we set the following
outflow boundary condition:
\[
\varepsilon \ds \big(\frac{1}{Re} \frac{\partial \u}{\partial \n} - p \n \big) = 0.
\]
where $\n$ denotes the outer normal. This last boundary condition results from the integration by parts when deriving the weak formulation, and it is called the do-nothing boundary condition. \\
%
%The iterative algorithm \eqref{Whi} is stopped when we have the following condition: for a tolerance $tol$,
%\[
%\ds \frac{|\u_h^{i+1} - \u_h^{i}|_X + \|p_h^{i+1} - p_h^i \|_{L^2(\Omega)}}{|\u_h^{i+1}|_X + \|p_h^{i+1} \|_{L^2(\Omega)}} \le tol.
%\]
%
%

In this section, we take $\f=0$ and we have
$$
\g = \left\{
\begin{array}{rcl}
\0 \quad && \mbox{ on } \Gamma_w \\
\u_{in} && \mbox{ on } \Gamma_{in},
\end{array}
\right.
$$
then, we can take ${\G} = (C_{in} (1-y)y,0)$ since we have $\G\in X$ and $\div(\varepsilon \G)=0$. Thus the iterative scheme is the following: Find $(\bar \w_h^{i+1}, \bar p_h^{i+1})\in Y_{0h} \times M_{h}$ such that for all $(\v_h, q_h) \in Y_{0h}\times M_{h}$\\
\begin{equation}\label{Whiu}
\left\{
\begin{array}{ll}
\medskip
\ds   a_h(\bar \w_h^{i+1}+{\G},\v_h)
+d_h(\tilde{\w}_h^i+{\G},\bar \w_h^{i+1}+{\G},\v_h) + (\beta (\varepsilon_h) | \bar{\w}_h^i+{\G}| (\bar{\w}_h^{i+1}+{\G}) , \v_h)
- b_h(\v_h,\bar p_h^{i+1}) = 0,\\
\ds b_h(\bar \w_h^{i+1}, q_h) = 0,
\end{array}
\right.
\end{equation}
where
$$
\tilde  \w^{i}_h = \ds \frac{ \bar \w_h^{i}+ \tilde \w^{i-1}_h}{2}
$$
and
$$
Y_{0h} =  X_{h} \cap \{ \w|_{\Gamma_{in} \cup \Gamma_w} = \0   \}.
$$
The corresponding indicators are:
\begin{equation}
\begin{array}{rcl}
\medskip
\eta_{i,\kappa}^L&=&   \|\bar \w_h^{i+1}- \bar \w_h^{i} \|_{H^1(\kappa)^d},\\
\medskip
\eta_{i,\kappa}^{D} &=& \ds h_{\kappa } || \frac{1}{Re}  \div (\varepsilon_h  \nabla  (\bar\w_h^{i+1}+{\G}_h)  ) - \alpha_{\varepsilon_h} (\bar \w_h^{i+1}+{\G}_h)  \ds - \varepsilon_h (  (\tilde\w_h^i+{\G}_h)\cdot \nabla )  ( \bar \w_h^{i+1}+{\G}) \\
\medskip%
&& \ds \qquad \qquad  -
\frac{1}{2} \div(\varepsilon_h (\tilde \w_h^i + {\G}_h) (\bar \w_h^{i+1} +{\G}_h)- \beta_{\varepsilon_h} |(\bar \w_h^i + +{\G}_h)| (\bar \w_h^{i+1}+{\G}_h)  - \varepsilon_h \nabla \bar p_h^{i+1} ||_{L^2(\kappa)^d}\\
\medskip
&& + \ds\frac{1}{2} \sum_{e \in \varepsilon_{\kappa}} h^{1/2}_{e} ||[(\frac{1}{Re} \varepsilon_h \nabla (\bar \w_h^{i+1} +{\G}_h)- \bar p_h^{i+1} \mathbb I)(\sigma)\cdot {\n}] ||_{L^2(e)^d}
 + ||\div (\varepsilon_h (\bar \w^{i+1}_h+{\G}_h))||_{L^2(\kappa)}.
\end{array}
\end{equation}
%
%In \cite{BDFtheornum1}, the authors treat a scheme similar to \eqref{Whiu} with $\bar \u_h^i$ instead of $\tilde \u_h^i$ and without the stabilization term $\ds \frac{1}{2} (\div (\bar \u_h^i) \w_h^{i+1}, \v_h)$. They study the dependency of the convergence of the iterative scheme \eqref{Whiu} with respect of $C_{in}$ for different values of the Reynolds number $Re$. They conclude that for a given tolerance $tol$, the bigger admissible value of $C_{in}$ denoted by $C_{in,max}$ decreases when $Re$ increases, and for all numerical test cases considered, the value of $C_{in,min}$ behaves like $C Re^{-1}$ ($C$ is a positive real constant).  \\

In the following, we will show numerical results corresponding to the iterative scheme \eqref{Whiu}. The main idea is to compare the uniform and adaptive methods. All the numerical results of this section are performed with an initial guess $\bar \u_h^0=\0$ and with a uniform initial mesh containing $4 N^2$ triangles with $N=60$. The algorithm is stopped by using the stopping criterion \eqref{eq:error_tolerance}.\\

\begin{figure}[htbp]
\begin{minipage}[b]{0.450\linewidth}
 \centering
\includegraphics[width=7cm]{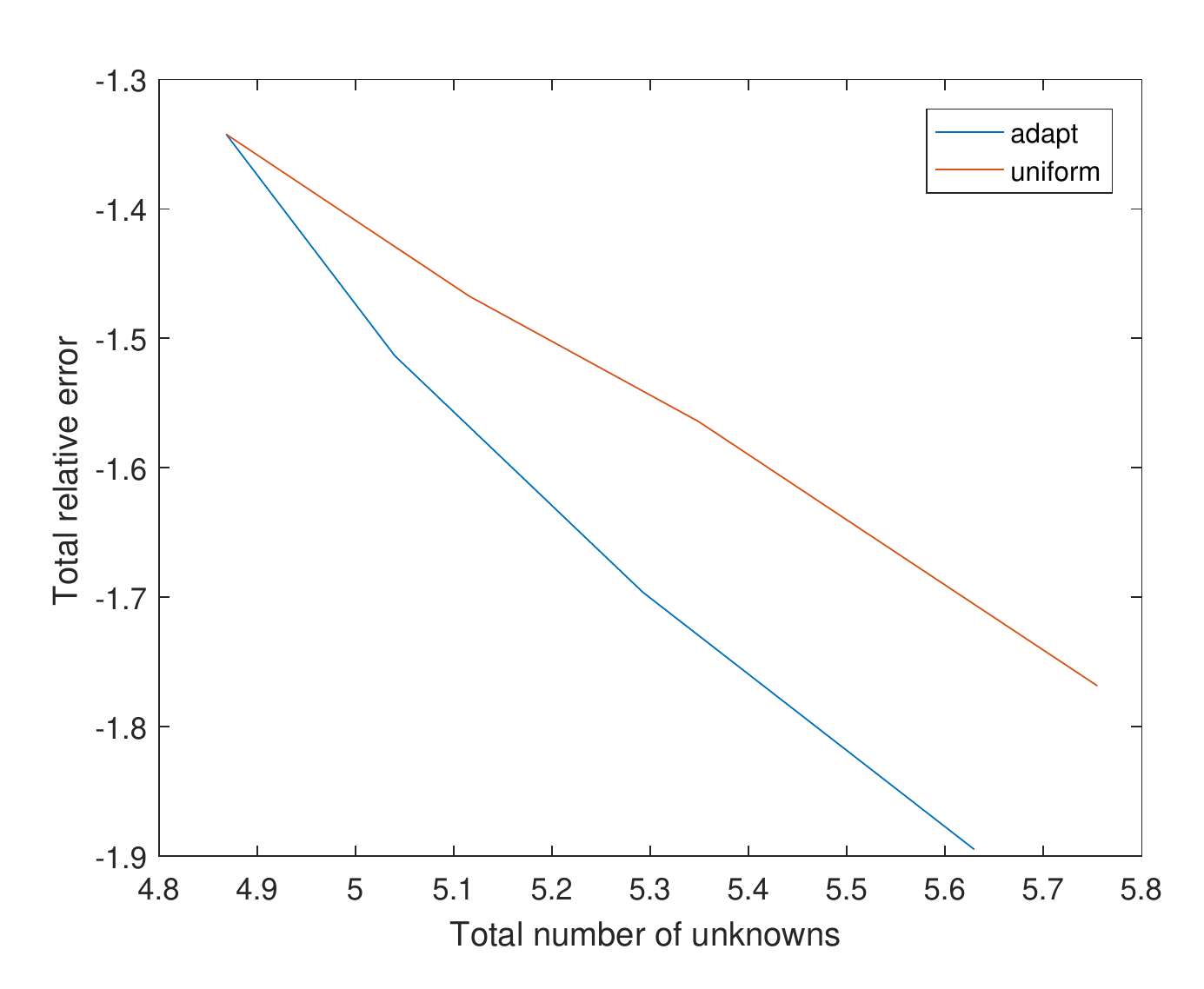}
\vspace{-.5cm}
\caption{$Re=100$, $C_{in}=.2$}
\label{figure11}
\end{minipage}%\hfill
\begin{minipage}[b]{0.450\linewidth}
\centering
\includegraphics[width=7cm]{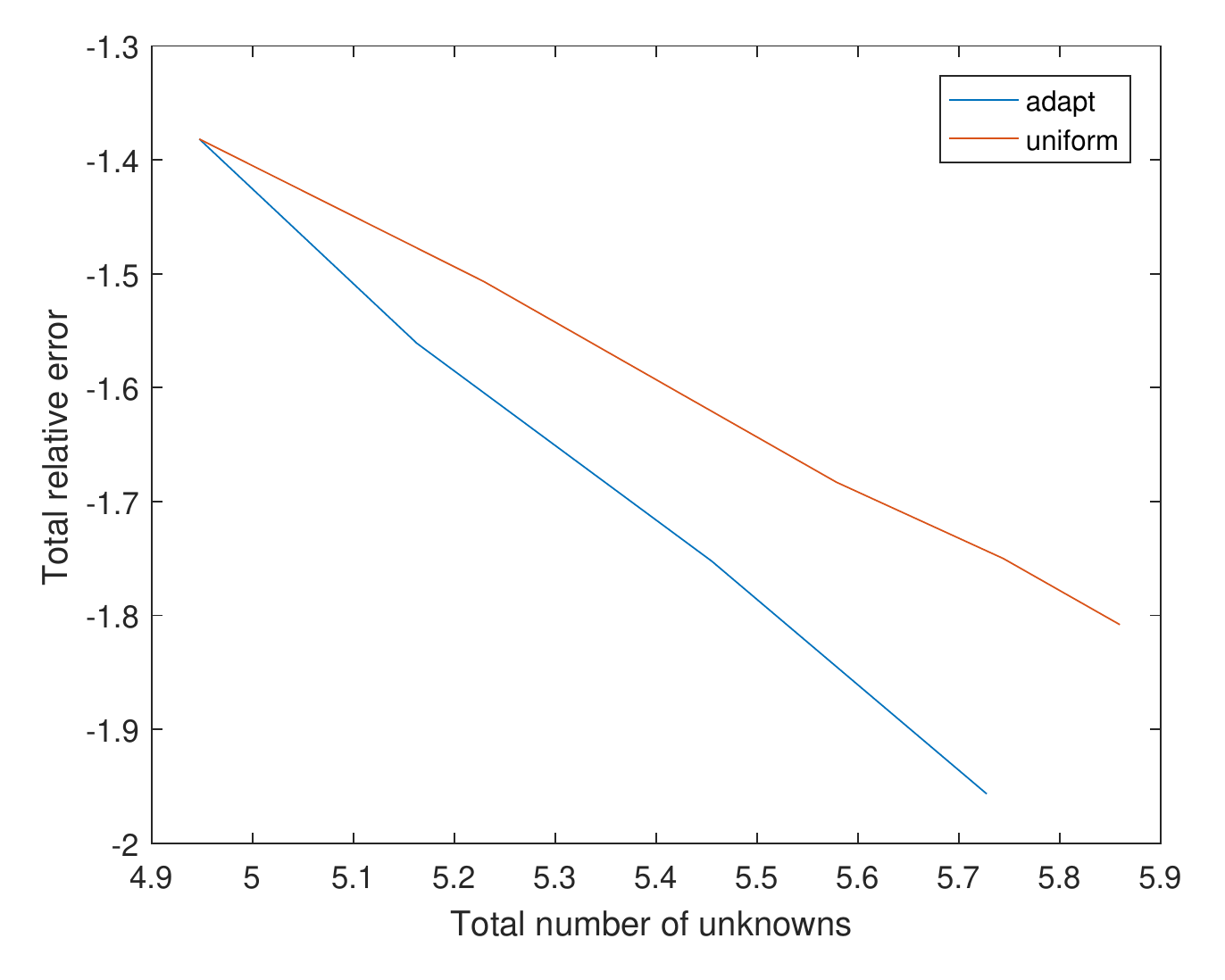}
\vspace{-.5cm}
\caption{$Re=500$, $C_{in}=.2$}
\label{ figure22}
\end{minipage}\hfill
\begin{minipage}[b]{0.450\linewidth}
\centering
\includegraphics[width=7cm]{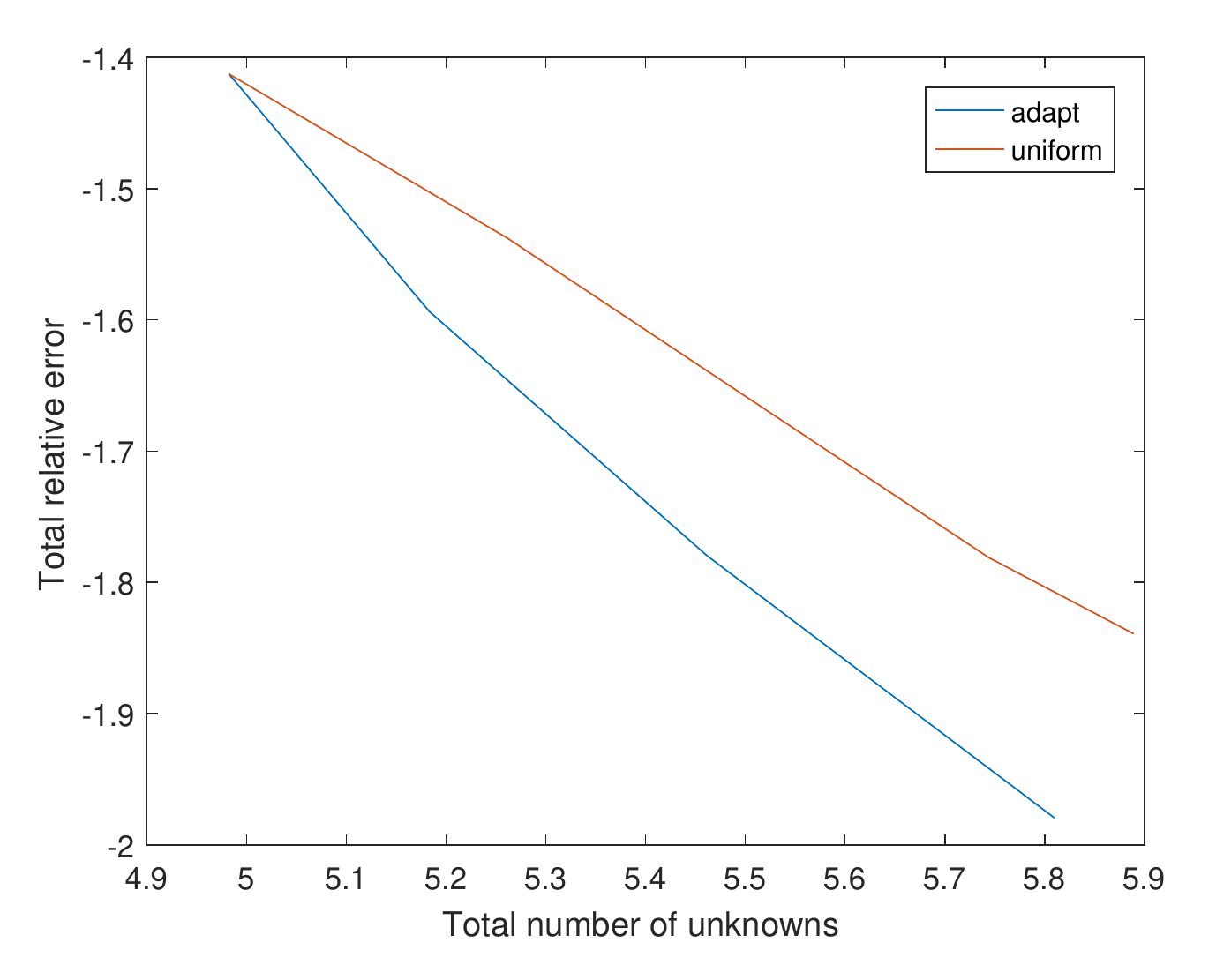}
\vspace{-.5cm}
\caption{$Re=1000$, $C_{in}=.2$} \label{figure33}
\end{minipage}%\hfill
\begin{minipage}[b]{0.450\linewidth}
\centering
\includegraphics[width=7cm]{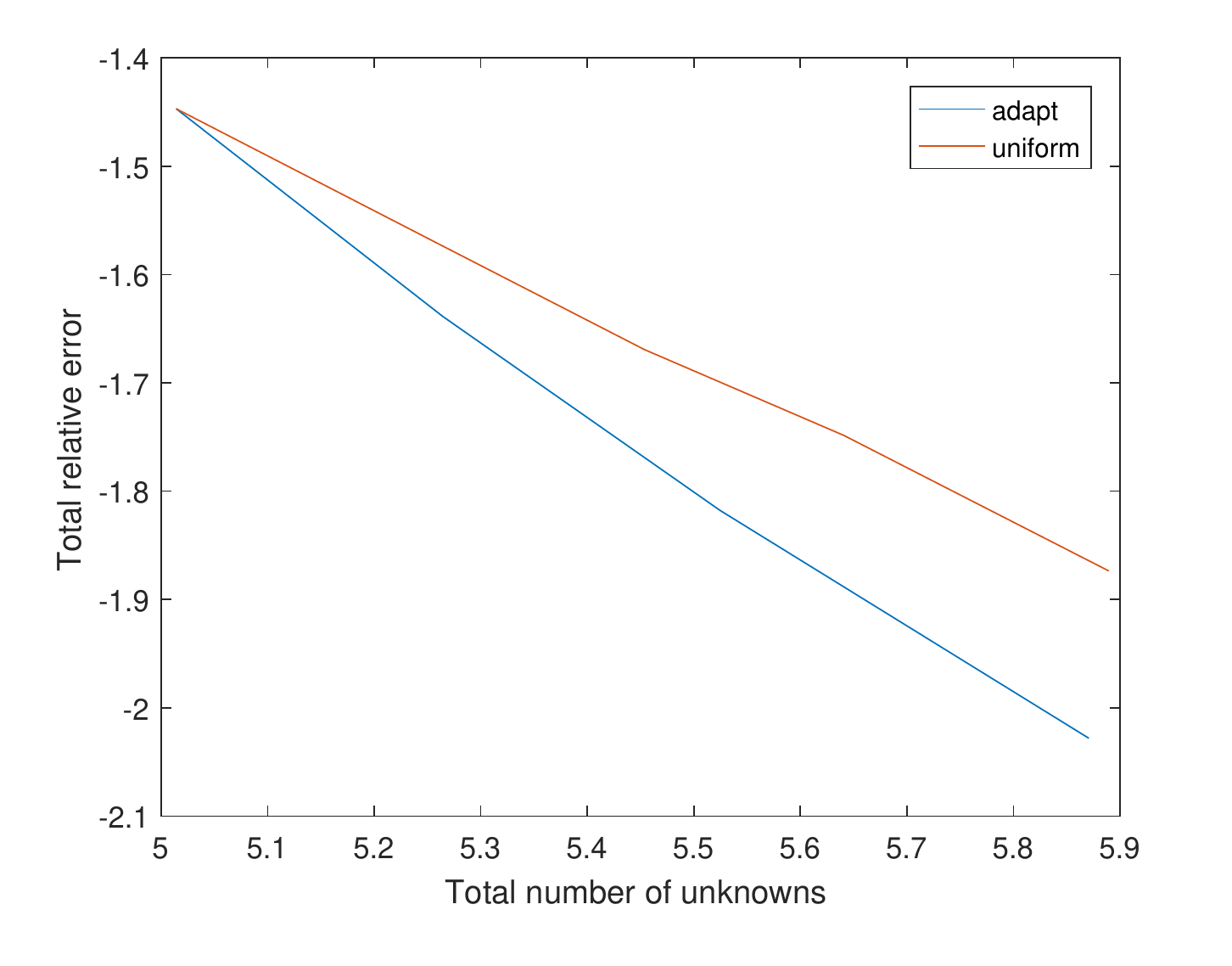}
\vspace{-.5cm}
\caption{$Re=2000$, $C_{in}=.2$} \label{figure44}
\end{minipage}\hfill
\end{figure}

Figure \ref{figure11}-\ref{figure44} show the total error indicator with respect of the total number of unknowns for $Re=100,500,1000,2000$ and for $C_{in}=0.2$. We can see clearly the advantage of the adapt method versus the uniform one.
\section{Conclusion} \label{sec:conclusion}
\noindent In this work, we have derived {\it {a posteriori}} error estimates for the finite element discretization of the Brinkman-Darcy-Forchheime system. These estimates yield an upper bound of the error which is computable up to unknown constants and allows to distinguish the discretization and the linearization errors. In this work, we show the advantages of the adaptive mesh refinement versus the uniform mesh method.

%\clearpage
\newpage
%%%%%%%%%%%%%%%%%%%%%%%%%%%%%%%%%%%%%%%%%%%%%%%%%%%%%%%%%%%%%%%%%%%%%%%%%
%%%%%%%%%%%%%%%%%%%%%%%%%%%%%%%%%%%%%%%%%%%%%%%%%%%%%%%%%%%%%%%%%%%%%%%%%
%%%%%%%%%%%%%%%%%%%%%%%%%%%%%%%%%%%%%%%%%%%%%%%%%%%%%%%%%%%%%%%%%%%%%%%%
%%%%%%%%%%%%%%%%%%%%%%%%%%%%%%%%%%%%%%%%%%%%%%%%%%%%%%%%%%%%%%%%%%%%%%%%
%

%
\end{document}